\documentclass[10pt]{amsart}
\usepackage{latexsym, amsmath, amssymb}
\usepackage{version}
\usepackage{color}
\usepackage[T1]{fontenc}
\usepackage{amsthm, amsopn, amsfonts}

\newif\ifpdf
\ifx\pdfoutput\undefined
   \pdffalse           
\else
   \pdfoutput=1
   \pdftrue            
\fi

\ifpdf
   \usepackage[pdftex]{graphicx}
\else
   
   \usepackage{graphicx}
\fi

\newcommand{\cal}{\mathcal}

\newcommand{\newsection}[1]
{\section{#1}\setcounter{theorem}{0} \setcounter{equation}{0}
\par\noindent}

\newtheorem{theorem}{Theorem}

\newtheorem{lemma}[theorem]{Lemma}

\newtheorem{proposition}[theorem]{Proposition}

\newcommand{\N}{{\mathbb N}}

\newcommand{\R}{{\mathbb R}}
\newcommand{\Z}{{\mathbb Z}}

\newcommand{\ang}{{\not\negmedspace\nabla}}

\newcommand{\e}{\epsilon}
\newcommand{\rs}{{r^*}}

\renewcommand{\S}{{\mathbb S}}
\newcommand{\M}{\cal M}
\newcommand{\tv}{{\tilde{v}}}
\newcommand{\tu}{{\tilde{u}}}
\newcommand{\tg}{{\tilde{g}}}
\newcommand{\tchi}{{\tilde{\chi}}}

\newcommand{\weight}{{\Bigl(1-\frac{2M}{r}\Bigr)}}

\begin{document}

\title
{
Strichartz estimates on Schwarzschild black hole backgrounds
}

\thanks{The authors were supported in part by the NSF 
grants DMS0354539 and DMS0301122
}

\author{Jeremy Marzuola}
\author{Jason Metcalfe}
\author{Daniel Tataru}
\author{Mihai Tohaneanu}

\address{Department of Applied Physics and Applied Mathematics,
  Columbia University, New York, NY 10027}

\address{Department of Mathematics, University of North Carolina,
  Chapel Hill, NC 27599-3250}

\address{Department of Mathematics, University of California,
  Berkeley, CA 94720-3840}

\address{Department of Mathematics, University of California,
  Berkeley, CA 94720-3840}

\begin{abstract}
  We study dispersive properties for the wave equation in the
  Schwarzschild space-time. The first result we obtain is a local
  energy estimate. This is then used, following the spirit of
  \cite{MT}, to establish global-in-time Strichartz
  estimates. A considerable part of the paper is devoted to a precise
  analysis of solutions near the  trapping region, namely the
   photon sphere.
\end{abstract}

\includeversion{1}
\includeversion{2}

\maketitle


\newsection{Introduction}

The aim of this article is to contribute to the understanding of the
global-in-time dispersive properties of solutions to wave equations on
Schwarzschild black hole backgrounds.  Precisely, we consider two
robust ways to measure dispersion, namely the local energy estimates
and the Strichartz estimates.

Let us begin with the local energy estimates. For solutions to the constant
coefficient wave equation in $3+1$ dimensions,
\[
\Box u = 0, \qquad u(0) = u_0, \quad u_t(0) = u_1,
\]
we have the original estimates of Morawetz
\cite{M}\footnote{There is another estimate commonly referred to as a
  Morawetz estimate.  This corresponds to using the multiplier
  $(t^2+r^2)\partial_t + 2tr\partial_r$.  We will reserve the term
  Morawetz estimate for \eqref{Morawetz} and shall call the latter
  estimate the Morawetz conformal estimate.}
\begin{equation}\label{Morawetz}
  \int_0^t\int_{\R^3} \frac{1}{|x|}|\ang u|^2(t,x)\:dt\:dx \lesssim
 \|\nabla u_0\|_{L^2}^2 + \| u_1\|_{L^2}^2
\end{equation}
where $\ang$ denotes the angular derivative.  To prove this one
multiplies the wave equation by the multiplier $(\partial_r+\frac{1}r) u$ and
integrates by parts.  Within dyadic spatial regions one can also
control $u$, $\partial_t u$ and $\partial_r u$.  Precisely, we have the
local energy estimates
\begin{equation}
R^{-\frac12} \|\nabla u\|_{L^2 (\R \times B(0,R))} 
+ R^{-\frac32} \| u\|_{L^2 (\R \times B(0,R))} 
\lesssim \|\nabla
u_0\|_{L^2} + \| u_1\|_{L^2}.
\label{localenergyflat}\end{equation}
See for instance \cite{KSS}, \cite{KPV}, \cite{SmSo}, \cite{St}, \cite{Strauss}.

One can also consider the inhomogeneous problem,
\begin{equation}
\Box u = f, \qquad u(0) = u_0, \quad u_t(0) = u_1,
\label{Minhom}\end{equation}
In view of \eqref{localenergyflat} we define the local energy space
$LE_M$ for the solution $u$ by
\begin{equation}
\|u\|_{LE_M} =  \sup_{j\in \Z} \Bigl[2^{-\frac{j}2}\|\nabla u\|_{L^2 (A_j)}
+ 2^{-\frac{3j}2}\|u\|_{L^2 (A_j)}\Bigr] ,
\label{leinf}\end{equation}
where
\[ 
A_j = \R\times \{2^j\le |x|\le 2^{j+1}\}. .
\]
For the inhomogeneous term $f$ we introduce a dual type norm
\[
\| f\|_{LE^*_M} = \sum_{j \in \Z} 2^{\frac{j}2} \|f\|_{L^2(A_j)}.
\]
Then we have:

\begin{theorem} 
The solution $u$ to \eqref{Minhom} satisfies the following estimate:
\begin{equation}
\| u\|_{LE_M} \lesssim  \|\nabla
u_0\|_{L^2} + \| u_1\|_{L^2} + \| f\|_{LE^*_M}
\label{Mest}\end{equation}
\end{theorem}

One may ask whether similar bounds  also hold for perturbations of the
Minkowski space-time.  Indeed, in the case of small long range
perturbations the same bounds as above were established very recently
by two of the authors, see  \cite[Proposition 2.2]{MT1} or
\cite[(2.23)]{MS2} (with no obstacle, $\Omega = \emptyset$). 
 See also \cite{Alinhac}, \cite{MS}
for related local energy estimates for small perturbations of the
d'Alembertian.  For large perturbations one faces additional
difficulties, due on one hand to trapping for large frequencies and on
the other hand to eigenvalues and resonances for low frequencies.
The Schwarzschild space-time, considered in the present paper, is a
very interesting example of a large perturbation of the Minkowski
space-time, where trapping causes significant difficulties.

The Schwarzschild space-time $\M$ is a spherically symmetric solution
to Einstein's equations with an additional Killing vector field $K$,
which models the exterior of a massive spherically symmetric body.
Factoring out the $\S^2$ component it can be represented via the
Penrose diagram:

\begin{figure}[h]
\centerline{ \scalebox{1}{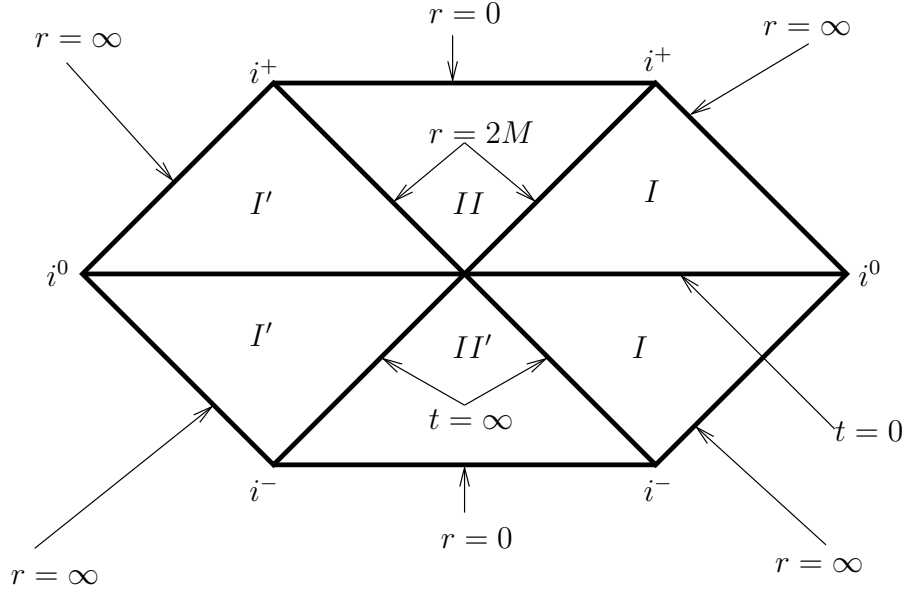}}
\caption{The Penrose diagram for the Kruskal extension
  of the Schwarzschild solution}
\label{f:pd}
\end{figure}

The radius $r$ of the $\S^2$ spheres is intrinsically determined and
is a smooth function on $\M$ which has a single critical point at the
center. The regions $I$ and $I'$ represent the exterior of the black
hole, respectively its symmetric twin, and are characterized by the
relation $r > 2M$. We can represent $I$ as
\[
I = \R\times (2M,\infty)\times \S^2
\]
with a metric whose line element is
\begin{equation}
ds^2=-\weight dt^2+\weight^{-1}dr^2+r^2d\omega^2
\label{swm}\end{equation}
where $d\omega^2$ is the measure on the sphere $\S^2$.  The
Killing vector field $K$ is given by $K = \partial_t$, which is
time-like within $I$. The differential $dt$ is intrinsic, but the function
$t$ is only defined up to translations on $I$. 

The regions $II$ and $II'$ represent the black hole, respectively its
symmetric twin, the white hole, and are characterized by the relation
$r < 2M$. The same metric as in \eqref{swm} can be used.  The Killing
vector field $K$ is still given by $K = \partial_t$, which is now
space-like.  Light rays can enter the black hole but not leave it. By
symmetry light rays can leave the white hole but not enter it.

The surface $r=2M$ is called the event horizon.  While the singularity
at $r=0$ is a true metric singularity, we note that the apparent
singularity at $r=2M$ is merely a coordinate singularity. Indeed,
denote
\[
r^*=r+2M\log(r-2M)-3M-2M\log M,
\]
so that
\[
d\rs = \weight^{-1} dr, \qquad \rs(3M)=0
\]
and set $v = t+r^*$. Then in the $(r,v,\omega)$ coordinates the metric
in region $I$ is expressed in the form
\[
ds^2=-\weight dv^2+2 dvdr+r^2d\omega^2,
\]
which extends analytically into the black hole region $I+II$.
In particular, given a choice of the function $t$ in region $I$, this
uniquely determines the function $t$ in the region $II$ via the
same change of coordinates.

In a symmetric fashion we set $w=t-\rs$.  Then in the $(r,w,\omega)$
coordinates the metric is expressed in the form
\[
ds^2=-\weight dw^2-2 dwdr+r^2d\omega^2,
\]
which extends analytically into the white hole region $I+II'$.

One can also introduce global nonsingular coordinates by rewriting the
metric in the Kruskal-Szekeres coordinate system,
\[
v' = e^{\frac{v}{4M}}, \qquad w' = -e^{-\frac{w}{4M}}.
\]
However, this is of less interest for our purposes
here.  Further information on the Schwarzschild space can be found in
a number of excellent texts.  We refer the interested reader to, e.g.,
\cite{HE}, \cite{MTW}, and \cite{Wald}.

 As far as the results in this paper are concerned, for large $r$ the Schwarzschild space-time can be viewed as
a small perturbation of the Minkowski space-time. The difficulties in
our analysis are caused by the dynamics for small $r$, where trapping
occurs.  The presence of trapped rays, i.e. rays which do not escape
either to infinity or to the singularity $r=0$, are known to be a
significant obstacle to proving local energy, dispersive, and
Strichartz estimates and, in some case, are known to necessitate a
loss of regularity.  See, e.g., \cite{BGT} and \cite{Ralston}.

There are two places where trapping occurs on the Schwarzschild
manifold.  The first is the surface $r=3M$ which is called the photon
sphere.  Null geodesics which are initially tangent to the photon
sphere will remain on the surface for all times.  Microlocally the
energy is preserved near such periodic orbits.  However what allows
for local energy estimates near the photon sphere is the fact that
these periodic orbits are hyperbolic.  The second is at the event
horizon $r=2M$, where the trapped geodesics are the vertical ones in
the $(r,v,\omega)$ coordinates.  However, this second family of
trapped rays turns out to cause no difficulty in the decay estimates
since in the high frequency limit the energy decays exponentially
along it as $v \to \infty$. This is due to the fact that the frequency
decays exponentially along the Hamilton flow, and in the physics
literature it is well-known as the red shift effect.

To describe the decay properties of solutions to the wave equation
in the Schwarzschild space, it is convenient to use coordinates which
make good use of the Killing vector field and are nonsingular along
the event horizon. The $(r,v,\omega)$ coordinates would satisfy these
requirements. However the level sets of $v$ are null surfaces,
which would cause some minor difficulties. This is why in $I + II$ we
introduce the function $\tilde v$ defined by
\[
\tv = v - \mu(r)
\]
where $\mu$ is a smooth function of $r$.
In the $(\tv,r,\omega)$ coordinates the metric has the form
\begin{multline*}
ds^2=-\weight d\tv^2 +2\left(1-\weight\mu'(r)\right) d\tv dr \\+
    \Bigl(2 \mu'(r) - \weight (\mu'(r))^2\Bigr)  dr^2  +r^2d\omega^2.
\end{multline*}

On the function $\mu$ we impose the following two conditions:

(i) $\mu (r) \geq  \rs$ for $r > 2M$, with equality for $r >
{5M}/2$.

(ii)  The surfaces $\tv = const$ are space-like, i.e.
\[
\mu'(r) > 0, \qquad 2 - \weight \mu'(r) > 0.
\]
\noindent The first condition (i) insures that the $(r,\tv,\omega)$
coordinates coincide with the $(r,t,\omega)$ coordinates in $r
>{5M}/2$. This is convenient but not required for any of our
results. What is important is that in these coordinates the metric
is asymptotically flat as $r \to \infty$. In the proof of the
Strichartz estimates, it is also required that $\mu'(r) =
\weight^{-1}$ near $r=3M$, which in other words says that we 
can work in the $(r,t)$ coordinates near the photon sphere.
However, this may be merely an artifact of our
method.

\begin{figure}[h]
\centerline{\scalebox{1}{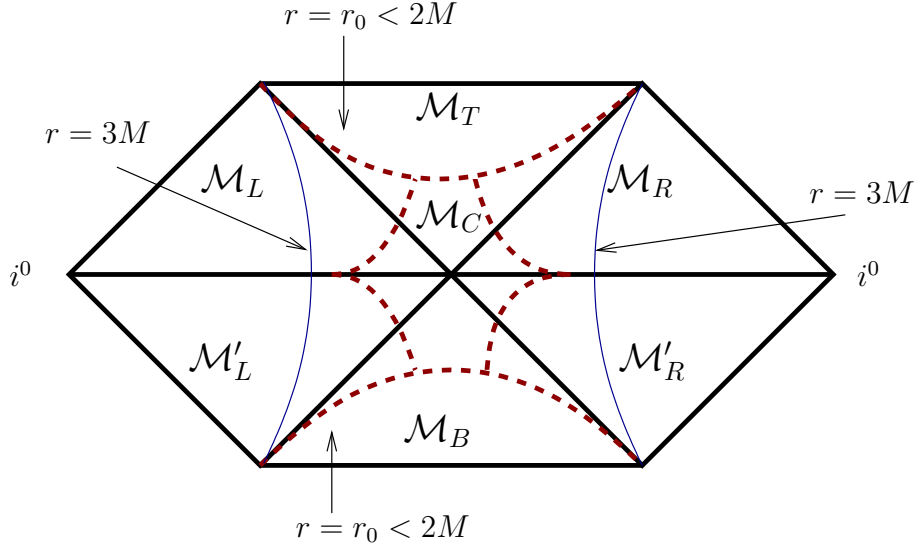}}
\caption{The Schwarzschild space partition represented on the Penrose diagram}
\label{f:pdcutoff}
\end{figure}

We introduce a symmetric function $\tilde{v}_1$ in $I'+II$, as well as
the the functions $\tilde{w}$ and $\tilde{w}_1$ in $I+II'$,
respectively $I'+II'$. Given a parameter $0 < r_0 < 2M$ we partition
the Schwarzschild space into seven regions
\[
\M = \M_R \cup \M_L \cup \M'_{R} \cup \M'_{L} \cup \M_T \cup \M_C
\cup \M_{B}
\]
as in Figure~\ref{f:pdcutoff}.
 The right/left top/bottom regions are
\[
\M_{R} = \{ \tv \geq 0, \ r \geq r_0 \} \subset I + II, \qquad \M_{L} = \{
\tv_1 \geq 0, \ r \geq r_0 \} \subset I' + II,
\]
\[
\M'_{R} = \{ \tilde w \leq 0, \ r \geq r_0 \} \subset I + II', \qquad \M'_{L} =
\{ \tilde w_1 \leq 0, \ r \geq r_0 \} \subset I' + II',
\]
the top and bottom regions are
\[
\M_T = \{ r < r_0\} \cap II, \qquad \M_B = \{ r < r_0\} \cap II',
\]
and the central region $\M_C$ is the remainder of $\M$. Moreover,
define
\[
\Sigma_R^- = \M_{R} \cap \{ \tv = 0 \}
\]
\[
\Sigma_R^+ = \M_R \cap \{ r = r_0\}.
\]
and similarly for the other regions.

In what follows we consider the Cauchy problem
\begin{equation}
\Box_g \phi = f, \qquad \phi_{|\Sigma_0} = \phi_0, \qquad \tilde K \phi_{|\Sigma_0} = \phi_1
\label{CP}\end{equation}
where for convenience we choose the initial surface $\Sigma_0$ to be
the horizontal surface of symmetry
\[
\Sigma_0 = \{ t = 0 \} \cap(I+I')
\]
\noindent
and $\tilde K$ is smooth, everywhere timelike and equals $K$ on
$\Sigma_0$ outside $\M_C$. Observe that we cannot use $K$ on all of 
$\Sigma_0$ since it is degenerate at the center (i.e. on the bifurcate sphere).

The equation \eqref{CP} can be solved as follows:

(i) Solve the equation in $\M_C$ with Cauchy data on $\Sigma_0$. Since
$\M_C$ is compact and has forward and backward space-like boundaries,
this is a purely local problem.

(ii) Solve the equation in $\M_R$ with Cauchy data on $\Sigma_R^-$.
The forward boundary of $\M_R$ is $\Sigma_R^+$, which is space-like.
This is the most interesting part, where we are interested in the
decay properties as $\tv \to \infty$. In a similar manner solve the
equation in $\M_L$, $\M'_R$ and $\M'_L$.

(iii) Solve the equation in $\M_T$ with initial data on the space-like
surface $\Sigma_T = \{ r=r_0\} \cap II$. Here one can track the solution
up to the singularity and encounter a mix of local and global
features. This part of the analysis in not pursued in the present article.

A significant role in our analysis is played by the Killing vector
field $K$, which in the $(r,\tv)$ coordinates equals $\partial_{\tv}$.  This is time-like outside the black hole but space-like
inside it. Furthermore, it is degenerate at the center.  Using the
Killing vector field outside the black hole we obtain a conserved
energy $E_0[\phi]$ for solutions $\phi$ to the homogeneous equation
$\Box_g \phi = 0$.  On surfaces $t=const$ in the $(r,t)$ coordinates
the energy $E_0[\phi](t)$ has the form
\begin{equation}\label{energy}
E_0[\phi] = \int_{S^2} \int_{2M}^\infty \Bigl[\weight^{-1}(\partial_t\phi)^2
+\weight (\partial_r \phi)^2+|\ang \phi|^2\Bigr]\ r^2 dr d\omega.
\end{equation}

Since the vector field $K$ is degenerate at the center, so is the
corresponding energy $E_0$ at $r=2M$. Hence it would be natural to replace
it with a nondegenerate energy, which on the initial surface
$\Sigma_0$ can be expressed as
\begin{equation}\label{Sigma_0}
E[\phi](\Sigma_0) = \int_{S^2} \int_{2M}^\infty
\Bigl[\weight^{-\frac32}(\partial_t\phi)^2 +\weight^\frac12
(\partial_r \phi)^2+|\ang \phi|^2\Bigr]\ r^2 dr d\omega.
\end{equation}
Unfortunately this is no longer conserved, and this is one of the
difficulties which we face in our analysis. We remark that a related form of a 
nondegenerate energy expression was
introduced in \cite{DR} and proved to be bounded in the exterior
region on surfaces $t = const$.

Part of the novelty of our approach is to prove bounds not only
in the exterior region, but also inside the event horizon. This is
natural if one considers the fact that the singularity 
at $r=2M$ is merely a removable coordinate singularity. In order 
to do this,  it is no longer suitable to measure the evolution of the 
energy on the surfaces $t=const$ (see below). Thus the above energy
$E[\phi](\Sigma_0)$ is relegated to a secondary role here and is used 
only to measure the size of the initial data.

A priori the energy $E[\phi](t)$ of $\phi$ only determines its Cauchy
data at time $t$ modulo constants. However, in what follows we
implicitly assume that $\phi$ decays at $\infty$, in which case
$\phi$ can be also estimated via a Hardy-type inequality,
\begin{equation}
 \int  \weight^{-\frac12} r^{-2} \phi^2\ r^2 dr d\omega
  \lesssim \int \weight^\frac12 (\partial_r \phi)^2 \ r^2 dr d\omega.
\label{hardy}\end{equation}
This is proved in a standard manner; the details are left to the reader.

We shall now further describe our main estimates in the region
$\M_{R}$: the local energy decay, the WKB analysis which yields a
local energy decay with only a logarithmic loss, and finally the
Strichartz estimates.

For the initial energy on $\Sigma_R^-$ we use
\[
E[\phi](\Sigma_R^-) = \int_{\Sigma_R^-}
 \left(
|\partial_r \phi|^2 +   |\partial_\tv \phi|^2    +
|\ang \phi|^2 \right) r^2  dr  d\omega.
\]
For the final energy on  $\Sigma_R^+$ we set
\[
E[\phi](\Sigma_R^+) = \int_{\Sigma_R^+}
 \left(  |\partial_r \phi|^2 +   |\partial_\tv \phi|^2    +
|\ang \phi|^2 \right) r_0^2   d\tv d\omega.
\]
We also track the energy on the space-like slices $\tv=const$,
\[
E[\phi](\tv_0) = \int_{ \M_R \cap \{\tv = \tv_0\}}
 \left(
|\partial_r \phi|^2 +   |\partial_\tv \phi|^2    +
|\ang \phi|^2 \right) r^2  dr  d\omega.
\]
Thus $E[\phi](\Sigma_R^-) = E[\phi](0)$.

For the local energy estimates one may first consider
a direct analogue of the Minkowski bound \eqref{Mest}. 
Unfortunately such a bound is hopeless due to the trapping which
occurs at $r=3M$.  Instead, for our first result we define a weaker
preliminary local energy space $LE_0$ with norm
\begin{equation} \label{le0}
\!\!\! \| \phi\|_{LE_0}^2  = \int_{\M_R} \left( \frac{1}{r^2}
|\partial_r  \phi|^2 +  \left(1-\frac{3M}r \right)^2 \Bigl(\frac{1}{r^2}|\partial_\tv \phi|^2    +
 \frac1r |\ang \phi|^2\Bigr) + \frac{1}{r^{4}} \phi^2 \right)  r^2  dr d\tv d\omega.
\end{equation}
Compared to the $LE_{M}$ norm we note the power loss in the angular
and $\tv$ derivatives at $r = 3M$. The $LE_0$ norm is also weaker
than $LE_M$ as $r \to \infty$, but this is merely for convenience.

At the same time we would like to  also consider the inhomogeneous
problem $\Box_g \phi = f$.  To measure the inhomogeneous term $f$, we
introduce the norm $LE_0^*$, which is stronger than $LE_M^*$:
\begin{equation} \label{le0stara}
\|f\|_{LE^*_{0}}^2 = \int_{\M_R}  \left(1-\frac{3M}r \right)^{-2} r^{2} f^2 \ r^2  dr d\tv d\omega.
\end{equation}
Again the important difference is at $r = 3M$. Our first local energy
estimate is the following:
\begin{theorem}\label{theorem.1}
  Let $\phi$ solve the inhomogeneous wave equation $\Box_g \phi = f$
  on the Schwarzs\-child manifold.  Then we have
\begin{equation}\label{main.estimate.inhom}
E[\phi](\Sigma_R^+) +\sup_{\tilde v \geq 0}
    E[\phi](\tilde v)+  \|\phi\|_{LE_0}^2 \lesssim E[\phi](\Sigma_R^-) +
\|f\|_{LE^*_0}^2.
\end{equation}
\end{theorem}
Here we made no effort to optimize the weights at $r=3M$ and
$r=\infty$. This is done later in the paper. On the other hand
the above estimate follows from a relatively simple application of the
classical positive commutator method. The advantage of having even
such a weaker estimate is that it is sufficient in order to allow
localization near the interesting regions $r = 3M$ and $r = \infty$,
which can then be studied in greater detail using specific tools.

 The first related results regarding the solution of the wave equation
 on Schwarzschild backgrounds were obtained in \cite{Wa} and \cite{KW}
 which proved uniform boundedness in region $I$ (including the event
 horizon). The first pointwise decay result (without, however, a rate
 of decay) was obtained in \cite{Tw}. Heuristics from 
 \cite{Pr} suggest that solutions to the wave equation in the
 Schwarzschild case should locally decay like $v^{-3}$. For spherically
 symmetric data a $v^{-3+\e}$ decay rate was obtained in \cite{DR4},
 and under the additional assumption of the initial data vanishing
 near the event horizon, the $v^{-3}$ decay rate was proved in
 \cite{Kr}. In general the best known decay rate, proved in \cite{DR},
 is $v^{-1}$ (see also \cite{BSt}). We also refer the reader
to \cite{SSS}, where optimal pointwise decay rates for each spherical harmonic
are established for a closely related problem.

Estimates related to \eqref{main.estimate.inhom} were first proved in
\cite{LS} for radially symmetric Schr\"odinger equations on
Schwarzschild backgrounds.  In \cite{BS1, BS2, BSerrata}, those
estimates are extended to allow for general data for the wave
equation.  The same authors, in \cite{BS3,BS4}, have provided studies
that give improved estimates near the photon sphere $r=3M$. 

Moreover, we note that variants of these bounds have played an important role in
the works \cite{BSt} and \cite{DR} which prove analogues of the
Morawetz conformal estimates on Schwarzschild backgrounds.  This
allows one to deduce a uniform decay rate for the local energy away
from the event horizon, though there is necessarily a loss of
regularity due to the trapping that occurs at the photon sphere.
Instead in this paper we restrict ourselves to time translation
invariant estimates, and we aim to clarify/streamline these as much as
possible.

All of the above articles use the conserved (degenerate) energy $E_0
[\phi]$ on time slices, obtained using the Killing vector field
$\partial_t$ . As such, their estimates are degenerate near the event
horizon. Further progress was made in \cite{DR}, where an additional
vector field was introduced near the event horizon, in connection to
the red shift effect. This enabled them to obtain bounds in the
exterior region involving a nondegenerate form of the energy related
to \eqref{Sigma_0}.

The approach of \cite{LS}, \cite{BS1}, \cite{BSt} and \cite{DR} is to
write the equation using the Regge-Wheeler tortoise coordinate and to
expand in spherical harmonics.  For the equation corresponding to each
spherical harmonic, one uses a multiplier which changes sign at
the critical point of the effective potential.

Here we work in the coordinates $(r,\tv,\omega)$, though this is not
of particular significance, and we do not expand into spherical
harmonics.  We prove \eqref{main.estimate.inhom} using a positive
commutator argument which requires a single differential multiplier.
 We hope that this makes the methods more robust for other
potential applications.

During final preparations of this article, localized energy estimates
proved without using the spherical harmonic decomposition also
appeared in \cite{DR2}.  The methods contained therein are somewhat
different from ours.

Compared to the stronger norms $LE_{M}$, $LE_{M}^*$ the weights in 
\eqref{main.estimate.inhom} have a polynomial singularity at $r = 3M$, 
which corresponds to the family of trapped geodesics on
the photon sphere.  As a consequence of the results we prove later,
see Theorem~\ref{theorem.3}, the latter fact can be remedied to
produce a stronger estimate.

\begin{theorem}\label{theorem.2}
  Let $\phi$ solve the inhomogeneous wave equation $\Box_g \phi = f$
  on the Schwarzs\-child manifold.  Then \eqref{main.estimate.inhom}
  still holds if the coefficient $(1-3M/r)^2$ in the $LE_0$ and the
  $LE^*_0$ norms is  replaced by
\[
 \left(1-\ln\left|1-\frac{3M}r\right|\right)^{-2}.
\]
\end{theorem}

Now we have only a logarithmic singularity at  $r=3M$.  The
result above is only stated in this form for the reader's convenience.
The full result in Theorem~\ref{theorem.3} is stronger but also more
complicated to state since it provides a more precise microlocal local
energy estimate.

The logarithmic loss is not surprising, since it is characteristic of
geometries with trapped hyperbolic orbits (see for instance
\cite{Burq}, \cite{Hans}, \cite{Zworski}). Indeed, a similar estimate
in the semiclassical setting is obtained in \cite{CVP} using entirely
different techniques. Note, however, that the aforementioned estimate
only involves logarithmic loss of the frequency; our result is
stronger since it also implies bounds for $\| (\ln |\rs|)^{-1} u\|_{L^2}$,
which are necessary in order to prove Strichartz estimates.

There are two regions on which the analysis is distinct.
The metric is asymptotically flat, and thus, near infinity,
one can retrieve the classical Morawetz type estimate.
On the other hand, around the photon sphere $r = 3M$ we take an
expansion into spherical harmonics as well as a time Fourier
transform. Then it remains to study an ordinary differential
equation which is essentially similar to
\[
(\partial_x^2 - \lambda^2(x^2+\epsilon)) u = f, \qquad |\epsilon|
\ll 1, \quad |x| \lesssim 1.
\]
For this we use a rough WKB approximation in the hyperbolic region
combined with energy estimates in the elliptic region. Airy type
dynamics occur near the zeroes of the potential.

Even though it is weaker, the initial bound in Theorem~\ref{theorem.1}
plays a key role in the analysis. Precisely, it allows us to glue
together the estimates in the two regions described above.

We next consider the Strichartz estimates.  For solutions to the
constant coefficient wave equation on $\R\times \R^3$, the
well-known Strichartz estimates state that
\begin{equation}\label{strich}
  \||D_x|^{-\rho_1} \nabla u\|_{L^{p_1}_tL^{q_1}_x} \lesssim
\|\nabla u(0)\|_{L^2}  + \||D_x|^{\rho_2} f\|_{L^{p_2'}_tL^{q_2'}_x}.
\end{equation}
Here the exponents $(\rho_i,p_i,q_i)$ are subject to the
scaling relation
\begin{equation}
\frac{1}p+\frac{3}q = \frac{3}2 -\rho
\label{scalingpq}\end{equation}
and the dispersion relation
\begin{equation}
\frac{1}p + \frac{1}q \leq \frac12, \qquad 2 < p \leq \infty.
\label{dispersionpq}\end{equation}
All pairs $(\rho,p,q)$ satisfying \eqref{scalingpq} and
\eqref{dispersionpq} are called Strichartz pairs. Those for which the
equality holds in \eqref{dispersionpq} are called sharp Strichartz
pairs.  Such estimates first appeared in the seminal works
\cite{Brenner}, \cite{Strichartz1, Strichartz2} and as stated include
contributions from, e.g., \cite{GV}, \cite{P}, \cite{K}, \cite{LiSo},
and \cite{KT}.

If one allows variable coefficients, such estimates are
well-understood locally-in-time.  For smooth coefficients, this was
first shown in \cite{MSS} and later for $C^2$ coefficients by
\cite{Smith} and \cite{T1,T2,T3}.

Globally-in-time, the problem is more delicate.  Even a small, smooth,
compactly supported perturbation of the flat metric may refocus a
group of rays and produce caustics. Thus, constructing a parametrix
for incoming rays proves to be quite difficult. At the same time, one
needs to contend with the possibility of trapped rays at high
frequencies and with eigenfunctions/resonances at low frequencies.

Global-in-time estimates were shown for small, long range
perturbations of the metric in \cite{MT} using an outgoing parametrix.
In order to keep the parametrix outgoing one must allow evolution both
forward and backward in time.  This construction is based on an
earlier argument in \cite{Tglobal} for the Schr\"odinger equation.
The smallness assumption, however, precludes trapping and does not
permit a direct application to the current setup.

On the other hand, a second result of \cite{MT} asserts that even for
large, long range perturbations of the metric one can still establish
global-in-time Strichartz estimates provided that a strong form of
the local energy estimates holds. This switches the burden to the
question of proving local energy estimates.

The result in \cite{MT} cannot be applied directly to the present
problem due to the logarithmic losses in the local energy estimates
near the trapped rays. However, it can be applied for the near
infinity part of the solution. In a bounded spatial region, on the
other hand, we take advantage of the local energy estimates to
localize the problem to bounded sets, in which estimates are shown
using the local-in-time Strichartz estimates of \cite{Smith},
\cite{T1}. Thus we obtain

\begin{theorem}
\label{Strichartz.theorem}
If $\phi$ solves $\Box_g \phi = f$ in $\M_R$ then for all nonsharp
Strichartz pairs $(\rho_1,p_1,q_1)$ and $(\rho_2,p_2,q_2)$ we have
\begin{equation}
 E[\phi](\Sigma_R^+) + \sup_{\tv} E[\phi](\tv) +  \left\| \nabla \phi\right\|^2_{L^{p_1}_\tv
\dot H^{-\rho_1,q_1}_x}
 \lesssim E[\phi](\Sigma_R^-) +
\left\| f \right\|^2_{L^{p_2'}_\tv \dot H^{\rho_2, q_2'}_x }.
\end{equation}
\end{theorem}
Here the Sobolev-type spaces $\dot H^{s,p}$ coincide with the usual
$\dot H^{s,p}$ homogeneous spaces in $\R^3$ expressed in polar
coordinates $(r,\omega)$.

As a corollary of this result one can consider the global solvability
question for the energy critical semilinear wave equation in the
Schwarzschild space,
\begin{equation}
\left\{ 
\begin{array}{lc} 
\Box_g \phi = \pm \phi^5  & \text{in } \M 
\cr \cr
\phi = \phi_0, \ \tilde K \phi = \phi_1 & \text{in } \Sigma_0.
\end{array}
\right.
\label{nonlin}\end{equation}

\begin{theorem}\label{tnlw}
Let $r_0 > 0$. Then there exists $\epsilon > 0$ so 
that for each initial data $(\phi_0,\phi_1)$ which satisfies
\[
E[\phi](\Sigma_0) \leq \epsilon
\]
the equation \eqref{nonlin} admits an unique solution
$\phi$ in the region $\{ r > r_0\}$ which satisfies the bound
\[
E[\phi](\Sigma_{r_0}) + \| \phi \|_{\dot H^{s,p}(\{r > r_0\})}  \lesssim
E[\phi](\Sigma_0)
\]
for all indices $s,p$ satisfying
\[
\frac{4}p = s+\frac12, \qquad 0 \leq s < \frac12.
\]
Furthermore, the solution has a Lipschitz dependence on the initial
data in the above topology. 
\end{theorem}

Some further clarification is needed for the function space 
$\dot H^{s,p}(\{r > r_0\})$ appearing above, in view of the ambiguity
due to the choice of coordinates. In a compact neighbourhood of
the center region $\M_C$  this is nothing but the classical 
$ H^{s,p}$ norm. By compactness, different choices of coordinates
lead to equivalent norms. Consider now the upper exterior region $\M_R$
(as well as its three other mirror images). Using the coordinates 
$(\tv,x)$ with $x =  \omega r$, we define $\dot H^{s,p}(\M_R)$ as the 
restrictions to $\R^+ \times \{|x| > r_0\}$ of functions in the 
homogeneous Sobolev space  $\dot H^{s,p}(\R \times \R^3)$.

{\bf Acknowledgements:} The authors are grateful to M. Dafermos and
I. Rodnianski for pointing out their novel way of taking advantage of
the red shift effect in \cite{DR}, and to N.  Burq and M. Zworski for
useful conversations concerning the analysis near trapped null
geodesics.


\newsection{The Morawetz-type estimate}

In this section, we shall prove Theorem \ref{theorem.1}. We note that
the estimate \eqref{main.estimate.inhom} is trivial over a finite
$\tv$ interval by energy estimates for the wave equation; the
difficulty consists in proving a global bound in $\tv$. By the same
token, once we prove \eqref{main.estimate.inhom} for some choice
of $r_0 < 2M$, we can trivially make the transition to any $r_0 < 2M$
due to the local theory. Thus in the arguments which follow we reserve
the right to take $r_0$ sufficiently close to $2M$.

We consider  solutions to the inhomogeneous wave equation on the
Schwarzschild manifold in $\M_R$, which is given by
\[
\Box_g\phi = \nabla^\alpha
\partial_\alpha \phi  = f.
\]
Here $\nabla$ represents the metric connection.  Associated to this
equation is an energy-momentum tensor given by
\[
Q_{\alpha\beta}[\phi]=\partial_\alpha\phi \partial_\beta \phi -
\frac{1}{2}g_{\alpha\beta}\partial^\gamma \phi \partial_\gamma \phi.
\]
A simple calculation yields the most important property of
$Q_{\alpha\beta}$, namely that if $\phi$ solves the homogeneous
wave equation then $Q_{\alpha\beta}[\phi]$  is divergence-free:
\[
\nabla^\alpha Q_{\alpha\beta}[\phi]=0,\quad \text{if }\nabla^\alpha
\partial_\alpha \phi =0.
\]
More generally, we have
\[
\nabla^\alpha Q_{\alpha\beta}[\phi]=\partial_\beta \phi \ \Box_g \phi.
\]

In order to prove Theorem \ref{theorem.1}, we shall contract
$Q_{\alpha\beta}$ with a vector field $X$ to form the momentum
density
\[
P_\alpha[\phi,X]=Q_{\alpha\beta}[\phi]X^\beta.
\]
Computing the divergence of this vector field, we have
\[
\nabla^\alpha P_\alpha[\phi,X] = \Box_g \phi X \phi+
Q_{\alpha\beta}[\phi]\pi^{\alpha\beta},
\]
where
\[
\pi_{\alpha\beta}=\frac{1}{2}(\nabla_\alpha X_\beta + \nabla_\beta X_\alpha)
\]
is the deformation tensor of $X$.

If $X$ is the Killing vector field $K$ then the above divergence
vanishes,
\begin{equation}
\nabla^\alpha P_\alpha[\phi,K] = 0 \quad \text{if}\quad \Box_g \phi = 0.
\label{killdiv}\end{equation}
This gives rise to the $E_0[\phi]$ conservation law outside the black
hole.

Naively, one may seek vector fields $X$ so that the quadratic form
$Q_{\alpha\beta}[\phi]\pi^{\alpha\beta}$ is positive definite.
However, this may not always be possible to achieve.  Instead we note
that it may be just as good to have the symbol of this quadratic form
positive on the characteristic set of $\Box_g$.  Then it would be
possible to make the above quadratic form positive after adding a
Lagrangian correction term of the form $ q
\partial^\gamma \phi\partial_\gamma \phi$. Such a term can be
conveniently expressed in divergence form modulo lower order terms.
Precisely, for a vector field $X$, a scalar function $q$ and a
$1$-form $m$ we define
\[
P_\alpha[\phi,X,q,m] = P_\alpha[\phi,X] + q \phi \partial_\alpha \phi
- \frac12 \partial_\alpha q \phi^2 + \frac{1}{2}m_{\alpha}\phi^2
\]
where $m$ allows us to modify the lower order terms in the divergence
formula. Then we obtain the modified divergence relation
\begin{equation}
\begin{split}
\nabla^\alpha P_\alpha[\phi,X,q,m] =&\  \Box_g \phi \Bigl(X\phi +
 q \phi\Bigr)+ Q[\phi,X,q,m],
\\ Q[\phi,X,q,m] = &\
Q_{\alpha\beta}[\phi]\pi^{\alpha\beta}
+ q  \partial^\alpha\phi\, \partial_\alpha \phi + m_\alpha
\phi\, \partial^\alpha\phi + \frac12(\nabla^\alpha m_\alpha -
\nabla^\alpha \partial_\alpha q) \, \phi^2.
\end{split}
\label{div}\end{equation}

Theorem~\ref{theorem.1} is proved by making appropriate choices for
$X$, $q$ and $m$ so that the quadratic form $Q[\phi,X,q,m]$ defined by
the divergence relation is positive definite. In what follows we
assume that $X$, $q$ and $m$ are all spherically symmetric and
invariant with respect to the Killing vector field $K$.

\begin{lemma}
  There exist smooth, spherically symmetric, $K$-invariant $X$, $q$,
  and $m$  in $r \geq 2M$ satisfying the following properties:

  (i) $X$ is bounded\footnote{In the $(r,\tilde v)$ coordinates} , $|q(r)| \lesssim r^{-1}$, $|q'(r)| \lesssim r^{-2}$
  and $m$ has compact support in $r$ .

(ii) The quadratic form $Q[\phi,X,q,m]$ is positive definite,
\[
Q[\phi,X,q,m] \gtrsim r^{-2} |\partial_r \phi|^2 + \left(1-\frac{3M}r
\right)^2 (r^{-2} |\partial_\tv \phi|^2 + r^{-1}|\ang \phi|^2) +
r^{-4} \phi^2.
\]

(iii) $X(2M)$ points toward the black hole, $X(dr)(2M) < 0$,
and $\langle m,dr\rangle(2M) > 0$.

\label{ibp}
\end{lemma}

We postpone the proof of the lemma and use it to conclude
the proof of Theorem~\ref{theorem.1}.  Let $X$, $q$ and $m$ be as in
the lemma.  We extend them smoothly beyond the event horizon
preserving the spherical symmetry and the $K$-invariance.
By \eqref{killdiv} we can modify the vector field $X$ without changing
the quadratic form $Q$ in \eqref{div},
\[
\nabla^\alpha P_\alpha[\phi,X + C K,q,m] =  \Box_g \phi \Bigl((X+CK)\phi +
 q \phi\Bigr)+ Q[\phi,X,q,m].
\]
Here $C$ is a large constant.
We integrate this relation in the region
\[
D = \{ 0 < \tv < \tv_0,\ r > r_0 \}
\]
using the $(r,\tv,\omega)$ coordinates. This yields
\begin{multline*}
\int_D \left( \Box_g \phi\, ((X+CK)\phi +
 q \phi)+ Q[\phi,X,q,m]\right) r^2 dr d\tv d\omega
\\=  \left. \int \langle d\tv, P[\phi,X + C K,q,m]\rangle
  r^2 dr  d\omega \right|_{\tv = 0}^{\tv = \tv_0}
  -  \int_{r=r_0} \langle dr,  P[\phi,X + C K,q,m]\rangle
  r_0^2 d\tv  d\omega.
\end{multline*}
We claim that if $C$ is large enough and $r_0$ sufficiently close to
$2M$ then the integrals on the right have the correct sign,
\begin{equation}
E[\phi](\tv_1) \lesssim -  \int_{\tv=\tv_1}
 \langle d\tv, P[\phi,X + C K,q,m]\rangle  r^2 dr  d\omega \lesssim
C E[\phi](\tv_1), \qquad v_1 \in \R
\label{tvsign}\end{equation}
\begin{equation}
  \langle d r, P[\phi,X + C K,q,m]\rangle  \gtrsim
|\partial_r \phi|^2 + |\partial_\tv \phi|^2 +  |\partial_\omega
\phi|^2 +\phi^2, \qquad r = r_0.
\label{tvsign2}\end{equation}
If these bounds hold then the conclusion of the theorem follows by
(ii) and Cauchy-Schwarz.

Indeed, a direct computation yields
\[
\langle d \tv,P[\phi,\partial_{\tv}] \rangle  = - \frac12 \left[\left( 2
    \mu' -\weight \mu'^2\right)
|\partial_{\tv}
\phi|^2 + \weight |\partial_r  \phi|^2  +
r^{-2} |\partial_\omega \phi|^2 \right],
\]
respectively
\[
\langle d r,P[\phi,\partial_{\tv}] \rangle =  |\partial_{\tv}\phi|^2 +
\weight (\partial_r - \mu' \partial_{\tv}) \phi  \partial_\tv \phi.
\]
On the other hand
\[
\langle d\tv,P[\phi,\partial_{r}] \rangle = \left(1-\weight \mu'\right)
|\partial_r \phi|^2  - \left(2 \mu'-\weight \mu'^2\right) \partial_\tv \phi
\partial_r \phi,
\]
while
\[
\langle dr,P[\phi,\partial_{r}] \rangle =  - \frac12 \left[- \left( 2
    \mu' -\weight \mu'^2\right)
|\partial_{\tv}
\phi|^2 - \weight |\partial_r  \phi|^2  +
r^{-2} |\partial_\omega \phi|^2 \right].
\]

We compute
\[
 \langle d\tv, P[\phi,X + C K]\rangle =
(X(d\tv)+ C) \langle d\tv,P[\phi,\partial_{\tv}] \rangle
+ X(dr) \langle d\tv,P[\phi,\partial_{r}] \rangle.
\]
For large enough $C$ we have $X(d\tv) + C \gtrsim C$.  Therefore the
first term on the right is negative definite for $r > 2M$. More
precisely, it is only the coefficient of the $|\partial_r \phi|^2$
term which degenerates at $r = 2M$. However, due to
condition (iii) in the lemma we have $X(dr)(2M) < 0$; therefore
we pick up a negative $|\partial_r \phi|^2$ coefficient at $r = 2M$.
Thus we obtain
\[
- \langle d\tv, P[\phi,X + C K]\rangle \approx  C \left[
|\partial_{\tv}
\phi|^2 + \weight |\partial_r  \phi|^2  +
r^{-2} |\partial_\omega \phi|^2 \right] +  |\partial_r  \phi|^2,
\qquad r > 2M.
\]
Since all the coefficients in the quadratic form on the left are
continuous, it follows that the above relation extends to $r > r_0$
for some $r_0 < 2M$ depending on $C$, namely
\begin{equation}
0 < 2M -r_0 \ll C^{-1}.
\label{rorange}\end{equation}

In order to prove \eqref{tvsign} it remains to estimate
the lower order terms $P[\phi,0,q,m]$ in terms of the
positive contribution above.  Since $|q| \lesssim r^{-1}$ and
$m$ has compact support in $r$, we can bound
\[
| \langle d\tv, P[\phi,0,q,m]\rangle| \lesssim  r^{-1} |\phi| \left[|\partial_{\tv}
\phi|^2 + |\partial_r  \phi|^2  
\right]^\frac12 + r^{-2}|\phi|^2.
\]
Then by Cauchy-Schwarz it suffices to estimate
\[
\int_{r_0}^\infty  r^{-2}|\phi|^2 r^2 dr  \lesssim    C^{-\frac12}
\int_{r_0}^\infty \left[ C \weight + 1\right]  |\partial_r \phi|^2 r^2 dr
\]
which is a routine Hardy-type inequality.

We next turn our attention to \eqref{tvsign2} and begin with the
principal part
\[
 \langle dr, P[\phi,X + C K]\rangle =
(X(d\tv)+ C) \langle dr,P[\phi,\partial_{\tv}] \rangle
+ X(dr) \langle dr,P[\phi,\partial_{r}] \rangle.
\]
Examining the expressions for the
two terms above, we see that for $r_0$ subject to \eqref{rorange}
we have
\[
 \langle dr, P[\phi,X + C K]\rangle \gtrsim C |\partial_\tv \phi|^2
+ |\partial_\omega \phi|^2  - \Bigl(1-\frac{2M}{r_0}\Bigr) |\partial_r \phi|^2, \qquad r = r_0.
\]
Next we consider the lower order terms. The contribution of
$m$ is
\[
\frac12 \langle m,dr \rangle \phi^2 \gtrsim \phi^2
\]
due to condition (iii) in the Lemma.
The contribution of $q$ is
\[
q \phi \langle dr, d\phi \rangle -\frac12 \phi^2  \langle dr, dq \rangle.
\]
The coefficient of the second term is $\weight q'$, which is negligible
for $r_0$ close to $2M$. In the first term
we have
\[
\langle dr, d\phi \rangle = \weight \partial_r \phi + \Bigl(1-\weight \mu'\Bigr)
\partial_\tv \phi.
\]
All terms involving $\weight$ are negligible, and since $q$ is bounded we get
\[
 q \phi \partial_\tv \phi \ll C |\partial_\tv \phi|^2 + \phi^2
\]
for large enough $C$.

\begin{proof}[Proof of Lemma~\ref{ibp}]
It is convenient to look for $X$ in the $(r,t)$ coordinates,
where we choose the vector field $X$ of the form
\[
X= X_1+ \delta X_2, \qquad \delta \ll 1
\]
where
\[
X_1 = a(r) \weight \partial_r, \qquad X_2 =
b(r) \weight \left( \partial_r - \weight^{-1} \partial_t\right)
\]
and $a$ and $b\weight$ will be chosen to be smooth.
Note that $X$ is a smooth vector field in the nonsingular
coordinates $(r,v)$, since in these coordinates we have
\[
X_1 = a(r)\left(\weight\partial_r+\partial_v\right), \qquad X_2 =
b(r)\weight\partial_r.
\]
We remark that the vector field $X_2$ is closely related to the vector field $Y$ introduced earlier
in \cite{DR} in order to take advantage of the red shift effect. However, in their construction $Y$ is in a form which is nonsmooth near the event horizon and which is restricted to the exterior region. 

The primary role played by $X_2$ here is to ensure that $X + CK$ is time-like near the event horizon.
The red-shift effect largely takes care of the rest.

For convenience, we set
\[
t_1(r) = \weight \frac{1}{r^2}\partial_r\Bigl(r^2 a(r)\Bigr).
\]
A direct computation yields
\begin{equation}
\label{div1}
\begin{split}
  \nabla^\alpha {P}_\alpha[\phi,X_1] = & \ \weight^2 a'(r) (\partial_r\phi)^2
+ a(r)\frac{r-3M}{r^2} |\ang \phi|^2 \\ & \ - \frac{1}{2}t_1(r)\partial^\gamma
\phi 
\,\partial_\gamma \phi+ X_1\phi \Box_g \phi,
\end{split}
\end{equation}
respectively
\begin{equation}
\label{div2}
\begin{split}
\nabla^\alpha {P}_\alpha[\phi,X_2] =  &\ \frac{1}{2}
b'(r)\left(\weight \partial_r \phi -\partial_t \phi\right )^2
\\ &\ + \left(  \frac{r-3M}{r^2} b(r) -\frac12 \weight b'(r)\right)
|\ang \phi|^2
\\ &\ - \frac{1}{r}\weight b(r)\partial^\gamma \phi\, \partial_\gamma
\phi + X_2\phi \Box_g \phi
\end{split}
\end{equation}
where
\[
\partial^\gamma \phi \, \partial_\gamma \phi =-
\left[\weight^{-1}(\partial_t \phi)^2 -\weight (\partial_r\phi)^2 -
  |\ang \phi|^2\right].
\]
We choose $a$ so that the first line of the right side of \eqref{div1}
is positive. This requires that
\begin{equation}
a'(r) \gtrsim r^{-2}, \qquad a(3M) = 0.
\label{areq}\end{equation}
We choose $b$ so that the first line of the right hand side of \eqref{div2}
is positive. Precisely, we take $b$ supported in $r \leq 3M$ with
\[
b = -\frac{b_0(r)}{1 -\frac{2M}r}, \qquad r \in [2M,3M].
\]
with $b_0$ smooth, decreasing in $[2M,3M)$ and supported in $\{r \leq
3M\}$. In particular that guarantees that $b_0(2M) > 0$, which is 
later used to verify the condition (iii) in the Lemma.

 The exact choice of $b_0$ is not important, and in effect $b$
only plays a role very close to the event horizon $r = 2M$.  Even
though $b$ is singular at $2M$, the second term of the coefficient of $|\ang \phi|^2$ in
the second line of \eqref{div2} is nonsingular. Hence if $ \delta$
is sufficiently small this term is controlled by the first line in
\eqref{div1}.

Taking the above choices into account, we have 
\begin{equation}
\begin{split}
Q[\phi,X,0,0] =&\ \weight^2 a'(r) (\partial_r\phi)^2 +
\delta\frac{1}{2} b'(r)\left(\weight
\partial_r \phi -\partial_t \phi \right)^2
\\ &\
+ O\left(\frac{(r-3M)^2}{r^3}\right)|\ang \phi|^2 - 
q_0\partial^\gamma \phi \, \partial_\gamma \phi
\end{split}
\label{firstq}\end{equation}
where 
\[
q_0(r) = \left(\frac{1}{2}t_1(r) + \delta\frac{1}{r}b(r)\weight \right).
\]

The last term in \eqref{firstq} is a Lagrangian expression and is
accounted for via the $q$ term. The first three terms give a
nonnegative quadratic form in $\nabla \phi$. This form is in effect positive 
definite for $r < 3M$, where $b' > 0$. However for larger $r$ it  
controls $\partial_r \phi$ and $\ang \phi$ but not $\partial_t \phi$.
This can be easily remedied with the Lagrangian term.
Precisely, we choose $q$ of the form
\[
q = q_0 +  \delta_1 q_1, \qquad q_1(r)= \chi_{\{ r > 5M/2\}} \frac{(r-3M)^2}{r^4}.
\]
where $\chi_{\{ r > 5M/2\}}$ is a smooth nonnegative cutoff 
which is supported in $\{ r > 5M/2\}$ and equals $1$ for $r > 3M$.
The positive parameter $\delta_1$ is chosen so that $\delta_1 \ll
\delta$. Then the only nonnegligible contribution of $\delta_1 q_1$ 
is the one involving $\partial_t \phi$. We obtain
\begin{equation}
\begin{split}
Q[\phi,X,q,0] =&\ \weight^2 O(r^{-2}) (\partial_r\phi)^2 +
\delta\frac{1}{2} 
b'(r)\left(\weight
\partial_r \phi -\partial_t \phi \right)^2
 \\ +& O\left(\frac{(r-3M)^2}{r^3}\right)|\ang \phi|^2 + \delta_1 q_1
\weight^{-1} |\partial_t \phi|^2  -   \frac12 \nabla^\alpha
\partial_\alpha q \, \phi^2.
\end{split}
\label{secondq}\end{equation}

The contribution of $q_1$ can be made arbitrarily small by taking
$\delta_1$ small. Hence it will be neglected in the sequel.
At this stage it would be convenient to be able to choose
$a$  so that $\nabla^\alpha \partial_\alpha   t_1(r) < 0$.
A direct computation yields
\[
\nabla^\alpha \partial_\alpha   t_1(r) = -L a
\]
with
\[
L a (r)= -\frac{1}{r^2}\partial_r\Bigl[\weight
r^2\partial_r\Bigl\{\weight \frac{1}{2r^2}
\partial_r\Bigl(r^2 a(r)\Bigr)\Bigr\}\Bigr].
\]
Unfortunately it turns out that the  condition $La > 0$ and \eqref{areq}
are incompatible, in the sense that there is no smooth
$a$ which satisfies both. However, one can find $a$ with
a logarithmic blow-up at $2M$ which satisfies both requirements.
Such an example is
\[
a(r) = r^{-2} \left( (r-3M)(r+2M) + 6 M^2\log\left(\frac{r-2M}M\right)\right).
\]
This is in no way unique, it is merely the simplest we were able
to produce. One verifies directly that
\[
a'(r) \gtrsim r^{-2}, \qquad L a (r) \gtrsim  r^{-4}.
\]
 To eliminate the singularity of $a$ above we replace it by
\[
a_\epsilon (r)= \frac{1}{r^2} f_{\epsilon}(R)
\]
where $\epsilon$ is a  small parameter,
\[
 R=(r-3M)(r+2M) + 6 M^2\log\left(\frac{r-2M}M\right),
\]
and
\[
f_{\epsilon}(R) = \epsilon^{-1} f(\epsilon R)
\]
where $f$ is a smooth nondecreasing function such that $f(R)=R$ on
$[-1, \infty]$ and $f = -2$ on $(- \infty,-3]$. The condition
\eqref{areq} is satisfied uniformly with respect to small $\e$;
therefore the choice of $\delta$ is independent of the choice
of $\epsilon$.

With this modification of $a$ we recompute
\[
L a_\epsilon = f'(\epsilon R) La + O(\epsilon) f''(\epsilon R)
+ O\left(\epsilon^2{\weight^{-1}} \right) f'''(\epsilon R).
\]
This is still positive except for the region $\{\e R < -1\}$.
To control it we introduce an $m$ term in the divergence
relation as follows.

Let $\gamma(r)$ be a function to be chosen later. We set
\[
m_t=\delta b'(r)\weight^2 \gamma, \qquad m_r= \delta b'(r)\weight
\gamma, \qquad m_\omega = 0.
\]
Then
\[
m_\alpha \partial^\alpha \phi = \delta b'(r) \gamma(r) \weight \left(\weight
\partial_r \phi -\partial_t \phi \right),
\]
while
\[
\nabla^\alpha m_\alpha = \delta r^{-2} \partial_r \left( \weight^2 r^2\,
  b'(r) \gamma (r)\right).
\]
Hence, completing the square we obtain
\[
\begin{split}
Q[\phi,X,q,m] =&\ \weight^2 O(r^{-2}) (\partial_r\phi)^2 +
O\left(\frac{(r-3M)^2}{r^3}\right)|\ang \phi|^2 \\
&\ + \delta_1 q_1 \weight^{-1} |\partial_t \phi|^2 + n \phi^2 \\
&\ + \delta\frac{1}{2} b'(r)\left(\weight
\partial_r \phi -\partial_t \phi +\weight \gamma \phi\right)^2 
\end{split}
\]
where the coefficient $n$ is given by
\[
\begin{split}
n = &\ L a_\epsilon - \frac{1}{2}\delta r^{-2} \partial_r r^2 \weight \partial_r
\left(r^{-1} b(r) \weight\right) - \delta \frac{b'(r)}2 \weight^2 \gamma(r)^2
\\ &\ + \frac{1}{2}\delta r^{-2} \gamma \partial_r \left(r^2 \weight^2
  b'(r)\right) +  \frac{1}{2}\delta \gamma'   \weight^2 b'(r).
\end{split}
\]
We assume that $\gamma$ is supported in $\{ r < 3M\}$ and satisfies
\[
0 \leq \gamma \leq 1, \qquad 
\gamma' > -1.
\]
Then for $r > 3M$ we have
\[
n =  L a_\epsilon \gtrsim r^{-4},
\]
while for $r \leq 3M$ we can write
\[
n =  L a_\epsilon +  \delta \gamma'(r)   \weight^2
  b'(r) + O(\delta).
\]
If $\epsilon R > -1$ then, using the bound from below on $\gamma'$,
 we further have
\[
n \geq L a + O(\delta)
\]
which is positive provided that $\delta$ is sufficiently small.
On the other hand in the region $\{\epsilon R \leq  -1\}$,
we have
\[
n \geq  \frac{1}{2}\delta 
 \weight^2
  b'(r)\gamma'(r)  + O(\delta) + O(\epsilon) f''(\epsilon R) +
 O\left(\epsilon^2{\weight^{-1}} \right) f'''(\epsilon R).
\]
The $\gamma'$ term can be taken positive, while all the other terms
may be negative so they must be controlled by it.  The restriction we
face in the choice of $\gamma'$ comes from the fact that $0 \leq
\gamma \leq 1$. Hence we need to verify that
\[
I = \int_{ \e R \leq -1} \delta + \epsilon |f''(\epsilon R)| +
\epsilon^2{\weight^{-1}}  |f'''(\epsilon R)| \ll \delta.
\]
Indeed, the interval of integration has size $\leq e^{-c
  \epsilon^{-1}}$; therefore the above integral can be bounded by
\[
I \lesssim e^{-c \epsilon^{-1}} + \epsilon
\]
which suffices provided that $\epsilon$ is small enough.

 Finally, note that
 \[
 X(dr)(2M)=\left(a(r) \weight + \delta b(r)
 \weight\right)(2M)<0,
 \]
 \[
 \langle m,dr\rangle(2M)=\left(\delta
 b'(r)\weight^2 \gamma\right)(2M)>0.
\]
So (iii) is also satisfied.
\end{proof}


\section{Log-loss local energy estimates}\label{secSchwlogLE}

The aim of this section is to prove a local energy estimate for
solutions to the wave equation on the Schwarzschild space which is
stronger than the one in Theorem~\ref{theorem.1}.  Consequently, we
strengthen the norm $LE_0$ to a norm $LE$ and we relax the norm
$LE^*_0$ to a norm $LE^*$ which satisfy the following natural bounds:
\begin{equation}
  \|  \phi \|_{LE_0}^2 \lesssim  \|  \phi \|_{LE}^2 \lesssim\|  \phi \|_{LE_{M}}^2,
\end{equation}
respectively
\begin{equation}
  \|  f \|_{LE_{M}^*}^2 \lesssim  \| f \|_{LE^*}^2 \lesssim\|  f \|_{LE_0^*}^2.
\end{equation}
We note that these bounds uniquely determine the topology of the $LE$
and $LE^*$ spaces away from the photon sphere and from infinity. This
is due to the fact that the local energy estimates in
Theorem~\ref{theorem.1} have no loss in any bounded region away from
the photon sphere. To define the $LE$, respectively $LE^*$, norms we 
consider a smooth partition of unity
\[
1 = \chi_{eh}(r) + \chi_{ps}(r) + \chi_\infty(r)
\]
where $\chi_{eh}$ is supported in $\{ r < 11 M/4\}$, $ \chi_{ps}$ is
supported in $\{ 5M/2 < r < 5M\}$ and $\chi_\infty$ is supported 
in $\{ r > 4M\}$. Then we set
\begin{equation}
\| \phi \|_{LE}^2 = \| \chi_{eh} \phi \|_{LE_{M}}^2 + \| \chi_{ps}
\phi \|_{LE_{ps}}^2 + \| \chi_{\infty} \phi \|_{LE_{M}}^2,
\end{equation}
respectively
\begin{equation}
\| \phi \|_{LE^*}^2 = \| \chi_{eh} \phi \|_{LE_{M}^*}^2 + \| \chi_{ps}
\phi \|_{LE^*_{ps}}^2 + \| \chi_{\infty} \phi \|_{LE^*_{M}}^2.
\end{equation}
The norms $LE_{ps}$ and $LE^*_{ps}$ near the photon sphere are defined
in Section \ref{pssec} below, see \eqref{leps}, respectively 
\eqref{leps*}; their topologies coincide with $LE_{M}$,
respectively $LE_{M}^*$, away from the photon sphere. 
 
With these notations, the main result of this section can be phrased 
in a manner similar to Theorem~\ref{theorem.1}:

\begin{theorem} \label{theorem.3} For all functions $\phi$ which solve
  $\Box_g \phi = f$ in $\M_R$ we have
\begin{equation}
\sup_{\tv > 0}  E[\phi](\tv) +  E[\phi](\Sigma_R^+) +
\| \phi \|_{LE}^2  \lesssim E[\phi](\Sigma_R^-) + \| f \|_{LE^*}^2.
\end{equation}
\end{theorem}

We continue with the setup and estimates near the photon sphere in
Section \ref{pssec}, the setup and estimates near infinity in Section \ref{inftysec}
and finally the proof of the theorem in Section \ref{tproof}.

\subsection{ The analysis near the photon sphere}
\label{pssec}

Here it is convenient to work in the Regge-Wheeler coordinates
given by
\[
r^*=r+2M\log(r-2M)-3M-2M\log M.
\]
Then $r = 3M$ corresponds to $\rs = 0$, and a neighbourhood of $r =
3M$ away from infinity and the event horizon corresponds to a compact
set in $\rs$.   In these coordinates  
the operator $\Box_g$ has the form
\begin{equation}
r \left (1-\frac{2M}r\right)  \Box_g r^{-1} = L_{RW} = \partial_t^2 - \partial_\rs^2 - \frac{r-2M}{r^3}
\partial_{\omega}  + V(r), \quad V(r) = r^{-1} \partial^2_\rs r.
\label{RW}\end{equation}

 For $\rs$ in a compact set the energy has the form
\[
E[\phi] \approx \int (\partial_t \phi)^2 + (\partial_\rs \phi)^2 +
(\partial_{\omega} \phi)^2 dr d \omega, 
\]
and the initial local smoothing norms are expressed as
\[
\|\phi\|_{LE_0}^2 \approx \int   (\partial_\rs \phi)^2  +
\rs^2 ( (\partial_{\omega} \phi)^2 +     (\partial_t \phi)^2 )   +
\phi^2  dr d \omega dt,
\]
respectively
\[
\|f\|_{LE_0^*}^2 \approx \int 
\rs^{-2} f^2  dr d \omega dt.
\]
On the other hand
\[
\|\phi\|_{LE_{M}}^2 \approx \int
  (\partial_\rs \phi)^2  +(\partial_{\omega} \phi)^2 +     (\partial_t
  \phi)^2 + \phi^2   dr d\omega dt,
\]
\[
\|f\|_{LE_{M}^*}^2 \approx \int  f^2  dr d \omega dt.
\]

In the sequel we work with spatial spherically symmetric pseudodifferential operators in the $(\rs,\omega)$ coordinates
where $\omega \in \S^2$. We denote by $\xi$ the dual variable 
to $\rs$, and by $\lambda$  the spectral parameter for $(-\partial_{\omega})^\frac12$. Thus the role of the Fourier variable
is played by the pair $(\xi,\lambda)$, and all our symbols are 
 of the form 
\[
 a(\rs,\xi,\lambda)
\]
To such a symbol we associate the corresponding Weyl operator $A^w$. 
Since there is no symbol dependence on $\omega$, one can view this
operator as a combination of a one dimensional Weyl operator and
the spectral projectors $\Pi_\lambda$ associated to the operator $(-\Delta_{\S^2})^\frac12$, namely
\[
 A^{w}= \sum_\lambda a^w(\lambda) \Pi_\lambda
\]
All of our $L^2$ estimates admit orthogonal decompositions with respect 
to spherical harmonics, therefore in order to prove them it suffices 
to work with the fixed $\lambda$ operators  $a^w(\lambda)$, and treat 
$\lambda$ as a parameter. However, in the proof of the Strichartz estimates
later on we need kernel bounds for operators of the form $A^w$,
which is why we think of $\lambda$ as a second Fourier variable
and track the symbol regularity with respect to $\lambda$ as well.
Of course, this is meaningless for $\lambda$ in a compact set; only 
the asymptotic behavior as $\lambda \to \infty$ is relevant.

Let $\gamma_0: \R \to \R^+$ be a smooth increasing function so that
\[
\gamma_0(y) = \left\{ \begin{array}{cc} 1& y < 1 , \cr y & y \geq 2.
  \end{array}
\right.
\]
Let $\gamma_1: \R^+ \to \R^+$ be a smooth increasing function so that
\[
\gamma_1(y) = \left\{ \begin{array}{cc} y^\frac12 & y < 1/2, \cr 1 & y \geq 1.
  \end{array}
\right.
\]
Let $\gamma: \R^2 \to \R^+$ be a smooth function with the following
properties:
\[
\gamma(y,z) = \left\{ \begin{array}{cc} 1 & z < C ,\cr 
\gamma_0(y)  &  y< \sqrt{z/2},\ z \geq C , \cr
z^\frac12 \gamma_1(y^2/z) & y \geq   \sqrt{z/2}, \   z   \geq C,
  \end{array}
\right.
\]
where $C$ is a large constant. In the sequel $z$ is a discrete
parameter, so the lack of smoothness at $z=C$ is of no consequence.

Consider the symbol
\[
a_{ps}(\rs,\xi,\lambda) = \gamma(-\ln(\rs^2 +\lambda^{-2} \xi^2), \ln \lambda),
\]
and its inverse 
\[
a_{ps}^{-1}(\rs,\xi,\lambda) = 
\frac{1}{\gamma(-\ln(\rs^2 +\lambda^{-2} \xi^2), \ln \lambda)}.
\]
We note that if $\lambda$ is small then they both equal $1$, while if
$\lambda$ is large then they satisfy the bounds
\begin{equation}
\begin{split}
  1  \leq a_{ps}(\rs,\xi,\lambda) \leq a_{ps}(\rs,0,\lambda)\leq (\ln
  \lambda)^\frac12,
  \\
  (\ln
  \lambda)^{-\frac12}\leq  a_{ps}^{-1}(\rs,0,\lambda) \leq a_{ps}^{-1}(\rs,\xi,\lambda) 
  \leq 1.
\end{split}
\label{apsbds}\end{equation}
We also observe that the region where $y^2 > z$ corresponds to $\rs^2
+\lambda^{-2} \xi^2 < e^{-\sqrt{ \ln \lambda}}$. Thus differentiating
the two symbols we obtain the following bounds 
\begin{equation}
|\partial_{\rs}^\alpha \partial_\xi^\beta \partial_\lambda^\nu
a_{ps}(\rs,\xi,\lambda)| \leq c_{\alpha,\beta,\nu} \lambda^{-\beta-\nu}
(\rs^2 +\lambda^{-2} \xi^2 + e^{-\sqrt{ \ln \lambda}} )^{-\frac{\alpha+\beta}2},
\label{apsbd}\end{equation}
respectively 
\begin{equation}
|\partial_{\rs}^\alpha \partial_\xi^\beta \partial_\lambda^\nu
a_{ps}^{-1}(\rs,\xi,\lambda)| \leq c_{\alpha,\beta,\nu}
a_{ps}^{-2}(\rs,\xi,\lambda) \lambda^{-\beta-\nu} (\rs^2 +\lambda^{-2}
\xi^2 + e^{-\sqrt{ \ln \lambda}} )^{-\frac{\alpha+\beta}2},
\label{aps1bd}\end{equation}
where $\alpha + \beta + \nu > 0$.  These show that we have a good
operator calculus for the corresponding pseudodifferential operators.
In particular in terms of the classical symbol classes we have
\[
a_{ps}, a_{ps}^{-1} \in                        
S^{\delta}_{1, 0}, \qquad \delta > 0.
\]

Then we introduce the Weyl operators
\[
A_{ps} = \sum_\lambda a_{ps}^w(\lambda) \Pi_\lambda,
\]
respectively
\[
A_{ps}^{-1} = \sum_\lambda (a_{ps}^{-1})^w(\lambda) \Pi_\lambda.
\]
By \eqref{apsbd} and \eqref{aps1bd} one easily sees that these 
operators are approximate inverses. More precisely for small
$\lambda$, $\ln \lambda < C$,
they are both the identity, while
for large $\lambda$ 
\[
\| a_{ps}^w(\lambda) (a_{ps}^{-1})^w(\lambda)-I\|_{L^2 \to L^2} \lesssim
\lambda^{-1} e^{\sqrt{\ln \lambda}}, \qquad \ln \lambda \geq C.
\]
Choosing $C$ large enough we insure that the bound above
is always much smaller than $1$.

We use these two operators in order to define the improved local
smoothing norms
\begin{equation}
  \| \phi\|_{LE_{ps}} =  \|A^{-1}_{ps} \phi\|_{H^1_{t,x}} \approx 
\| A_{ps}^{-1} \nabla_{t,x} \phi\|_{L^2},
\label{leps}\end{equation}
\begin{equation}
\| f \|_{LE^*_{ps}} = \| A_{ps} f\|_{L^2}.
\label{leps*}\end{equation}

Due to the inequalities \eqref{apsbds} we have a bound 
from above for $a_{ps}^w(\lambda)$,
\[
\|a_{ps}^w(\lambda) f\|_{L^2} \lesssim \| a_{ps} (\rs,0,\lambda)
f\|_{L^2} \lesssim \| |\ln |\rs| |  f\|_{L^2}, 
\]
respectively a bound from below for $(a_{ps}^{-1})^w(\lambda)$,
\[
\|(a_{ps}^{-1})^w(\lambda) f\|_{L^2} \gtrsim  \| a_{ps}^{-1} (\rs,0,\lambda)
f\|_{L^2} \gtrsim \| |\ln |\rs| |^{-1}  f\|_{L^2} 
\]
for $f$ supported near $\rs = 0$. In particular this shows
that for $f$  supported near the photon sphere we have
\begin{equation}
 \| \phi\|_{LE_{ps}} \gtrsim \| |\ln |\rs| |^{-1}  \nabla \phi
 \|_{L^2}, \qquad \| f \|_{LE^*_{ps}} \lesssim \| |\ln |\rs| |f \|_{L^2}
\end{equation}
which makes Theorem~\ref{theorem.2} a direct consequence of 
Theorem~\ref{theorem.3}.

Our main estimate near the photon sphere is

\begin{proposition}
a) Let $\phi$ be a function supported in $\{ 5M/2 < r < 5M \}$
which solves $\Box_g \phi = f$. Then
\begin{equation}
 \| \phi\|_{LE_{ps}}^2\lesssim \| f\|^2_{ LE_{ps}^*}.
\label{psa}\end{equation}

b) Let $f \in LE^*_{ps}$ be supported in $\{ 11M/4 < r < 4M \}$. Then there
is a function $\phi$ supported in $\{ 5M/2 < r < 5M \}$ so that
\begin{equation}
\sup_t E[\phi] + \| \phi\|_{LE_{ps}}^2 
+ \| \Box_g \phi -f\|_{LE_0^*}^2
\lesssim \| f\|^2_{ LE^*_{ps}}.
\label{psb} \end{equation}
\label{pps}\end{proposition}

\begin{proof}
  Due to \eqref{RW} we can recast the problem in the Regge-Wheeler
  coordinates. Denoting $u=r\phi$, $g = \weight r f$ we have $L_{RW}
  u = g$. Also it is easy to verify that for $\phi$ and $f$
  supported in a fixed compact set in $\rs$ we have
\[
\| \phi \|_{LE_{ps}} \approx \|u\|_{LE_{ps}}, \qquad \| f \|_{LE_{ps}^*} \approx \|g\|_{LE_{ps}^*}.
\]
Hence in the proposition we can replace $\phi$ and $f$ by $u$ and $g$,
and $\Box_g$ by $L_{RW}$.

To prove part (a) we expand in spherical harmonics with respect to the
 angular variable and take a time Fourier transform. We are left with
 the ordinary differential equation
\begin{equation}
(\partial_\rs^2 + V_{\lambda,\tau}(\rs)) u = g ,
\label{pseq}\end{equation}
where
\[
V_{\lambda,\tau}(\rs) = \tau^2 - \frac{r-2M}{r^3} \lambda^2 + V.
\]

Depending on the relative sizes of $\lambda$ and $\tau$ we
consider several cases. In the easier cases it suffices to replace the
bound \eqref{psa} with a simpler bound
\begin{equation}
  \| \partial_{\rs} u \|_{L^2} + (|\tau|+|\lambda|)\|u\|_{L^2} \lesssim
  \|g\|_{L^2}.
\label{psaeasy}\end{equation}

{\bf Case 1:} $\lambda,\tau \lesssim 1$.
Then we solve \eqref{pseq} as a Cauchy problem with data on one
side and obtain a pointwise bound ,
\[
|u|+|u_{\rs}| \lesssim \|g\|_{L^2}
\]
which easily implies \eqref{psaeasy}.

{\bf Case 2:} $\lambda \ll \tau$.  Then $V_{\lambda,\tau}(\rs) \approx
\tau^2$ for $\rs$ in a compact set; therefore \eqref{pseq} is
hyperbolic in nature.  Hence we can solve \eqref{pseq} as a Cauchy
problem with data on one side and obtain
\[
\tau |u|+|u_{\rs}| \lesssim \|g\|_{L^2} ,
\]
which implies \eqref{psaeasy}.

{\bf Case 3:} $\lambda \gg \tau$. Then $V_{\lambda,\tau}(\rs) \approx
-\lambda^2$ for $\rs$ in a compact set; therefore \eqref{pseq} is
elliptic.  Then we solve \eqref{pseq} as an elliptic problem with
Dirichlet boundary conditions on a compact interval  and obtain
\[
\lambda^{\frac32} |u|+\lambda^{\frac12} |u_{\rs}| \lesssim \|g\|_{L^2} ,
\]
which again gives \eqref{psaeasy}.

{\bf Case 4:} $\lambda \approx \tau \gg 1$.  In this case
\eqref{psaeasy} is no longer true, and we need to prove \eqref{psa}
which in this case can be written in the form
\begin{equation}
\label{psahard}
\| \partial_{\rs} u \|_{L^2} + \lambda \| (a_{ps}^{-1})^w(\lambda)  u\|_{L^2}
\lesssim \| a_{ps}^w(\lambda)  g\|_{L^2} ,
\end{equation}
where $u,g$ are subject to \eqref{pseq}.  The $\partial_{\rs} u$ term
above is present in order to estimate the high frequencies $|\xi| \gg
\lambda$. For lower frequencies it is controlled by the second term 
on the left of \eqref{psahard}.

The potential $V$ in \eqref{pseq} can be treated perturbatively in
\eqref{psa} and is negligible.  The remaining part of
$V_{\lambda,\tau}(\rs)$ has a nondegenerate minimum at $r = 3M$ which
corresponds to $\rs = 0$.  Hence we express it in the form
\[
V_{\lambda,\tau}(\rs) = \lambda^2(W(\rs) +\e) ,
\]
where $W$ is smooth and has a nondegenerate zero minimum at $\rs=0$
and $|\e| \lesssim 1$.

 We now prove the following:
\begin{proposition}\label{ODE}
Let $W$ be a smooth function satisfying  $W(0)=W'(0)=0$, $W''(0)>0$,
and $|\e| \lesssim 1$. Let $w$ be a solution of the ordinary differential equation
\[
\partial_\rs^2 + \lambda^2(W(\rs) +\e)) w(\rs) = g
\]
supported near $\rs = 0$. Then \eqref{psahard} holds.
\end{proposition} 

It would be convenient to replace the norm on the right 
in \eqref{psahard} by $ \| a_{ps}(\rs,0,\lambda)  g\|_{L^2}$. 
This is not entirely possible since this is a stronger norm. However,
we can split $g$ into a component $g_1$ with $a_{ps}(\rs,0,\lambda)
g_1 \in L^2$  plus a high 
frequency part:

\begin{lemma}
Each function $g \in L^2$ supported near the photon sphere
can be expressed in the form
\[
g = g_1 + \lambda^{-2} \partial_{\rs}^2 g_2
\]
with $g_1$ and $g_2$ supported near the photon sphere so that 
\begin{equation}
  \| a_{ps}(\rs,0,\lambda)  g_1\|_{L^2} + 
\||\rs^2 + e^{-\sqrt{\ln \lambda}}|^{\frac18} g_2\|_{L^2} +
 \lambda^{-2} \|\partial_{\rs} g_2\|_{L^2}
\lesssim    \| a_{ps}^w(\lambda)  g\|_{L^2}.
\label{gud}\end{equation}
\end{lemma}
\begin{proof}
The symbols $a_{ps}(\rs,0,\lambda)$ and $a_{ps}(\rs,\xi,\lambda)$
are comparable provided that 
\[
\ln (\rs^2 + e^{-\sqrt{\ln \lambda}})  \approx \ln (\rs^2 +
e^{-\sqrt{\ln \lambda}}+\lambda^{-2} \xi^2).
\]
This includes a region of the form 
\[
D = \left\{  \ln (\lambda^{-2} \xi^2) < 
\frac{1}{8} \ln (\rs^2 + e^{-\sqrt{\ln \lambda}}) \right\}.
\]
We note that the factor $\frac18$, arising also in the exponent of the
second term in \eqref{gud}, is somewhat arbitrary. A small choice
leads to a better bound in \eqref{gud}.

If $\chi$ is a smooth function which is $1$ in $(-\infty,-1]$ and
$0$ in $[0,\infty)$ then we  define a smooth 
characteristic function $\chi_D$ of the domain $D$ by
\[
\chi_D(\rs,\xi,\lambda) = \chi ( \ln (\lambda^{-2} \xi^2)  - \frac18 \ln (\rs^2 + e^{-\sqrt{\ln \lambda}})).
\]
One can directly compute the regularity of $\chi_D$,
\[
\chi_D \in S^{0}_{1,\delta}, \qquad \delta > 0.
\]
To obtain the decomposition of $g$ we set
\[
g_2 =   q^w g ,
\]
where the symbol of $q$ is 
\[
q(\rs,\xi,\lambda) = \lambda^2 \xi^{-2} (1-\chi_D).
\]
Since $(a_{ps}^{-1})^w$ is an approximate inverse for $a_{ps}^w$,
the estimate for $g_2$ in the lemma can be written in the form
\begin{equation}
\|(\rs^2 + e^{-\sqrt{\ln \lambda}})^{\frac18} q^w (a_{ps}^{-1})^w(\lambda) f\|_{L^2} +
 \lambda^{-2} \|\partial_{\rs}  q^w (a_{ps}^{-1})^w(\lambda) f\|_{L^2}
\lesssim    \| f\|_{L^2}.
\label{guda}\end{equation}
In the first term it suffices to look at the principal symbol of the
operator product since the remainder belongs to
$OPS^{-1+\delta}_{1,\delta}$ for all $\delta > 0$.  To verify that the
product of the symbols is bounded we note that $a_{ps}^{-1}$ is
bounded. For the other two factors we consider two cases.  If $|\xi|
\gtrsim \lambda$ then both factors are bounded. On the other hand
if $|\xi| \lesssim \lambda$ then in the support of $q$ we have
\[
\lambda^{-2} \xi^2 \gtrsim (\rs^2 + e^{-\sqrt{\ln \lambda}})^{\frac18} ,
\]
which gives 
\[
q \lesssim (\rs^2 + e^{-\sqrt{\ln \lambda}})^{-\frac18}.
\]
The estimate for the second term in \eqref{guda} is similar but
simpler.

It remains to consider the bound for $g_1$, which is given by 
\[
g_1 = (1 + \lambda^{-2} D_{\rs}^2 q^w) g, \qquad D_{\rs} = \frac{1}i 
\partial_\rs
\]
As above, the bound for $g_1$ can be written in the form
\[
\| a_{ps}(\rs,0,\lambda)  (1 + \lambda^{-2} D_{\rs}^2 q^w)
(a_{ps}^{-1})^w (\lambda) f\|_{L^2} \lesssim \|f\|_{L^2}.
\]
The three operators above belong respectively to $S^{\delta}_{1,\delta}$,
$S^0_{1,\delta}$, and $S^\delta_{1,\delta}$ for all $\delta > 0$. Hence 
the product belongs to $S^{\delta}_{1,\delta}$, and  it suffices to
show that its principal symbol is bounded. But the principal symbol of
the product is given by 
\[
a_{ps}(\rs,0,\lambda) \chi_D a_{ps}^{-1}(\rs,\xi,\lambda) ,
\]
which is bounded due to the choice of $D$.

Finally we remark that as constructed the functions $g_1$ and $g_2$
are not necessarily supported near the photon sphere. This is easily rectified by replacing them with truncated versions,
\[
g_1:= \chi_1(\rs) g_1, \qquad g_2:= \chi(\rs) g_2  
\]
where $\chi_1$ is a smooth compactly supported cutoff which equals $1$ in the support of $g$. It is clear that the bound \eqref{gud} is still valid after truncation.

\end{proof}

Using the above decomposition of $g$ we write $u$ in the form
\[
u = \lambda^{-2} g_2 + \tilde u.
\]
For the first term we use the above lemma to estimate
\[
\lambda \| \lambda^{-2} g_2 \|_{L^2} + \| \lambda^{-2} \partial_\rs g_2 \|_{L^2}
\lesssim \| a_{ps}^w (\lambda) g\|_{L^2} ,
\]
which is stronger than what we need. 
For  $\tu$ we write the equation
\begin{equation}
(\partial_\rs^2 + \lambda^2(W+\e)) \tilde{u} = \tilde g,
\qquad \tilde g = g_1 - (W+\e) g_2.
\end{equation}
For $\tg$ we only use a weighted $L^1$ bound,
\[
\| (\lambda^{-1}+|W+\e|)^{-\frac14} \tg\|_{L^1} \lesssim
\| a_{ps}^w(\lambda) g\|_{L^2},
\]
which is obtained from the weighted $L^2$ bounds on $g_1$ and $g_2$ 
by Cauchy-Schwarz.

For $\tu$ on the other hand, it suffices to obtain a pointwise bound:

\begin{lemma}
For each $\lambda^{-1} < \sigma < 1$ and each function
$\tu$ with compact support, we have
\[
\lambda \| (a_{ps}^{-1})^w(\lambda) \tu\|_{L^2} \lesssim \|
(\sigma+|W+\e|)^{-\frac14} \partial_{\rs} \tu\|_{L^\infty} + 
\lambda \|
(\sigma+|W+\e|)^{\frac14}  \tu\|_{L^\infty}.
\]
\end{lemma}
\begin{proof}
  Since $W$ has a nondegenerate zero minimum at $0$, if $\e > -\sigma$
  then $\sigma+|W+\e| \approx \sigma+|\epsilon|+ W$. Hence without any restriction  in generality we can replace $(\epsilon,\sigma)$ by
$(0, \sigma +|\epsilon|)$. Thus in the sequel we can assume that either 
 $\e = 0$ or $\e < -\sigma$. We
  consider three cases:

{\bf Case 1:} $|\e|,\sigma < e^{-\sqrt{\ln \lambda}} $. 
We consider an almost orthogonal partition of $\tu$ in dyadic regions
  with respect to $\rs$:
\[
\tu = \tu_{< s_0} + \sum_{s_0 \leq s <  1} \tu_s, \qquad 
s_0 = e^{-\frac12 \sqrt{\ln \lambda}}. 
\]
For each piece we can freeze  $\rs$ in the symbol 
of $a_{ps}^{-1}$ and  estimate in $L^2$
\[
\| (a_{ps}^{-1})^w (\lambda)\tu_s \|_{L^2} \approx \| a_{ps}^{-1}(s,D,\lambda)
\tu_s\|_{L^2}, \qquad \| (a_{ps}^{-1})^w (\lambda)\tu_{<s_0} \|_{L^2} \approx \| a_{ps}^{-1}(0,D,\lambda)
 \tu_{<s_0} \|_{L^2}.
\]

In the first case we use the symbol bound
\[
\lambda a_{ps}^{-1}(s,\xi,\lambda) \lesssim \lambda |\ln s|^{-1} + s^{-\delta} |\xi|,
\qquad \delta > 0 ,
\]
where the second term accounts for the region where $|\xi| > \lambda
s^\delta$. This yields
\[
\lambda \| (a_{ps}^{-1})^w(\lambda) \tu_s \|_{L^2} \lesssim  \
|\ln s|^{-1}  \lambda \|  \tu_s\|_{L^2} + s^{-\delta}
\|\partial_{\rs} \tu_s\|_{L^2}.  
\]
In the support of $\tu_s$ we have $\sigma+|W+\e| \approx s^2$; 
therefore by Cauchy-Schwarz we obtain
\[
\lambda \| (a_{ps}^{-1})^w (\lambda)\tu_s \|_{L^2} \lesssim  |\ln s|^{-1} 
\lambda\| (\sigma+|W+\e|)^{\frac14}  \tu\|_{L^\infty} +
s^{1-\delta}
\|(\sigma+|W+\e|)^{-\frac14} \partial_{\rs} \tu\|_{L^\infty}. 
\]
The summation with respect to $s$ follows due to  the $L^2$ almost orthogonality
of the functions $(a_{ps}^{-1})^w (\lambda)\tu_s$,
\[
\begin{split}
\lambda \Bigl\| (a_{ps}^{-1})^w (\lambda)\sum_{s_0\le s < 1}\tu_s \Bigr\|_{L^2}
 \lesssim  &\ \lambda \left( \sum_{s_0\le s < 1} \| (a_{ps}^{-1})^w(\lambda) \tu_s \|_{L^2}^2 \right)^\frac12
\\
\lesssim &\ 
 \lambda \|
(\sigma+|W+\e|)^{\frac14}  \tu\|_{L^\infty} + \|
(\sigma+|W+\e|)^{-\frac14} \partial_{\rs} \tu\|_{L^\infty} .  
\end{split}
\]
This orthogonality is due to the fact that the kernel
of $(a_{ps}^{-1})^w (\lambda)$ decays rapidly on the $\lambda^{-1+\delta}$
scale in $\rs$, therefore the overlapping of the two functions
of the form $(a_{ps}^{-1})^w(\lambda) \tu_s$ is trivially small for 
nonconsecutive values of $s$.

On the other hand for the center piece $ \tu_{< s_0}$ a similar
computation yields
\[
\begin{split}
\lambda \| (a_{ps}^{-1})^w (\lambda)\tu_{<s_0} \|_{L^2} \lesssim &\
|\ln s_0|^{-1}  \lambda \|  \tu_{<s_0}\|_{L^2} + s_0^{-\delta}
\|\partial_{\rs} \tu_{< s_0}\|_{L^2}  
\\ 
\lesssim
& \ \lambda 
\|(\sigma+|W+\e|)^{\frac14} \tu\|_{L^\infty} +  s_0^{1-\delta} \|
(\sigma+|W+\e|)^{-\frac14}  \partial_{\rs}  \tu\|_{L^\infty} ,
\end{split}
\]
where the weaker bound in the first term is due to the fact that 
\[
\int_{|\rs| < s_0} (\sigma+|W+\e|)^{-\frac12} d\rs \lesssim \ln \lambda.
\]

{\bf Case 2:}  $\e=0$, $\sigma \geq e^{-\sqrt{\ln \lambda}}$.
Then we select $s_0$ by $s_0^2 = \sigma$ and partition $\tu$
into 
\[
\tu = \tu_{< s_0} + \sum_{s_0 \leq s \leq  1} \tu_s.
\]
The analysis proceeds as in the first case, with the simplification
that there is no longer a singularity in the weight
$(\sigma+|W+\e|)^{-\frac12}$ for $|\rs| < s_0$.

{\bf Case 3:} $\sigma, e^{-\sqrt{\ln \lambda}} < -\e$.
Then we select $s_0$ by $s_0^2 = -\e$ and partition $\tu$
into 
\[
\tu = \tu_{< s_0} + \tu_{s_0} + \sum_{s_0 < s \leq  1} \tu_s.
\]
Then all pieces are estimated as in Case 2 with the exception 
of $ \tu_{s_0}$, where we have to contend with the singularity
in the weight. However, compared to Case 1 we have a better 
integral bound
\[
\int_{|\rs| \approx s_0} |W+\e|^{-\frac12} d\rs \lesssim 1
\]
and the conclusion follows again.
\end{proof}

Due to the above lemma, it suffices to prove that, given $| \e|
\lesssim 1$ and $(u,g)$ supported near the photon sphere so that
\begin{equation}
(\partial_\rs^2 + \lambda^2(W+\e)) u = g,
\label{pseqe} \end{equation}
then for some $\lambda^{-1} < \sigma < 1$ we have
\begin{equation}
 \|(\sigma+|W+\e|)^{-\frac14} \partial_{\rs} u\|_{L^\infty} + 
\lambda \|(\sigma+|W+\e|)^{\frac14}  u\|_{L^\infty}
\lesssim \| (\lambda^{-1}+|W+\e|)^{-\frac14} g\|_{L^1}. 
\label{pseqest}\end{equation} 
We remark that the first term on the left gives the $L^2$ bound for
$u_\rs$ in \eqref{psahard}.

We consider three subcases depending on the choice of $\epsilon$:

{\bf Case 4 (i)} $ \epsilon \gg \lambda^{-1}$. Then
$|W+\e| \approx \rs^2 +\e$. Choosing $\sigma = \epsilon$, it suffices to prove that:

\begin{lemma}
  Suppose that $\epsilon \gg \lambda^{-1}$. Then for $u$ with compact
  support solving \eqref{pseqe}, we have
\[
 \lambda (\epsilon+\rs^2)^\frac14  |u| +  (\epsilon+\rs^2)^{-\frac14}
|u_\rs| \lesssim \| (\e+\rs^2)^{-\frac14} g\|_{L^1}. 
\]
\end{lemma}

\begin{proof}
  We solve \eqref{pseqe} as a Cauchy problem from both sides toward
  $0$.  For this we use an energy functional which is inspired by the
  classical WKB approximation,
\begin{eqnarray*}
  E(u(\rs)) = \lambda^2 ( W+\epsilon)^{\frac{1}{2}} u^2 +
(W+\epsilon)^{-\frac{1}{2}} u_\rs^2 + \frac12
W_\rs ( W+\epsilon)^{-\frac{3}{2}} u u_\rs.
\end{eqnarray*}
Since $|W_\rs| \lesssim W^\frac12$, the condition $\e \gg \lambda^{-1}$
guarantees that $E$ is positive definite.  Computing its
derivative, we have
\[
\begin{split}
\frac{d}{d\rs} E(u(\rs)) = & \frac12 (W_{\rs\rs} ( W+\epsilon)^{-\frac{3}{2}} -\frac32 W_\rs^2
 ( W+\epsilon)^{-\frac{5}{2}})  u u_\rs
\\ & +  (2(W+\epsilon)^{-\frac{1}{2}} u_\rs + \frac12 W_\rs (
W+\epsilon)^{-\frac{3}{2}} u) g.
\end{split}
\]
This leads to the bound
\[
\left|\frac{d}{d\rs} E(u(\rs))\right| \lesssim  \lambda^{-1} (W+\epsilon)^{-\frac{3}{2}}
 E(u(\rs)) + E^\frac12(u(\rs))  (W+\e)^{-\frac14} g.
\]
Since
\[
\int \lambda^{-1} (W+\epsilon)^{-\frac{3}{2}} dr \lesssim 1,
\]
the conclusion follows by  Gronwall's inequality.
\end{proof}

{\bf Case 4 (ii)} $|\e| \lesssim \lambda^{-1}$. We choose $\sigma = \lambda^{-1}$. Then
$ \sigma + |W+\epsilon| \approx \lambda^{-1}+\rs^2$.
Hence it suffices to prove the pointwise bound:

\begin{lemma}
  Suppose that $|\epsilon| \lesssim \lambda^{-1}$. Then for $u$ with
  compact support solving \eqref{pseqe}, we have
\[
\lambda (\lambda^{-1}+\rs^2)^\frac14 |u| +
(\lambda^{-1}+\rs^2)^{-\frac14}
|u_\rs| \lesssim \| (\lambda^{-1}+\rs^2)^{-\frac14} g\|_{L^1} .
\]
\end{lemma}

\begin{proof}
  We use the same energy functional $E$ as above. However, this time
  $E$ is only positive definite when $W \gg \lambda^{-1}$ or
  equivalently $|\rs| \gg \lambda^{-\frac12}$. Applying Gronwall as in
  the first case yields the conclusion of the Lemma for $|\rs| \gg
  \lambda^{-\frac12}$.

  In the remaining interval $\{ |\rs| \lesssim \lambda^{-\frac12}\}$ we
  view \eqref{pseqe} as a small perturbation of the equation
  $\partial_\rs^2 u = g$.  Precisely, we can use the energy functional
\[
E_1(u(\rs)) = \lambda^\frac32 |u|^2 + \lambda^{\frac12} |u_\rs|^2 ,
\]
which satisfies
\[
\Bigl|\frac{d}{d\rs} E_1(u(\rs))\Bigr| \lesssim \lambda^{\frac12} E_1(u(\rs))
+   E_1^\frac12(u(\rs)) \lambda^{\frac14} g.
\]
This allows us to use Gronwall's inequality to estimate the
remaining part of $u$.
\end{proof}

{\bf Case 4 (iii)} $-\epsilon \gg \lambda^{-1}$. Then we 
choose $\sigma = |\epsilon|^{\frac13}\lambda^{-\frac23}$ and
prove the
pointwise bound:
\begin{lemma}
  Suppose that $-\epsilon \gg \lambda^{-1}$. Then for $u$ with compact
  support solving \eqref{pseqe} we have
\[
\lambda (|W+\epsilon| + |\epsilon|^{\frac13} \lambda^{-\frac23})^\frac14
|u| + (|W+\epsilon| + |\epsilon|^{\frac13}
\lambda^{-\frac23})^{-\frac14} |u_\rs| \lesssim \|(|W+\epsilon| +
|\epsilon|^{\frac13} \lambda^{-\frac23})^{-\frac14}g\|_{L^1}.
\]
\end{lemma}

\begin{proof}
The energy $E$ above is still useful for as long as it stays positive
definite, i.e. if
\begin{equation}
|W_\rs| < 2 \lambda |W + \epsilon|^{3/2}.
\label{wwrs}\end{equation}
Given the quadratic behavior of $W$ at $0$, this amounts to
\[
W+\epsilon \gg \lambda^{-\frac23} |\e|^{\frac13}.
\]
In this range, due to \eqref{wwrs}  we obtain, as in Case 4(i), 
\[
\Bigl|\frac{d}{d\rs} E(u(\rs))\Bigr| \lesssim  \lambda^{-1} (W+ \epsilon)^{-\frac{3}{2}}
 E(u(\rs)) + E^\frac12(u(\rs))  (W+ \e)^{-\frac14} g,
\]
which by Gronwall's inequality and Cauchy-Schwarz yields the desired
bound.

In a symmetric region around the zeroes of $W+\e$,
\[
|W+\e| \lesssim \lambda^{-\frac23} |\e|^{\frac13},
\]
the bounds for $u$ and $u_\rs$ remain unchanged and the equation
\eqref{pseqe} behaves like a small perturbation of $\partial_\rs^2 u =
g$, and we can use a straightforward modification of the argument above.

Finally, in the region
\[
[r_-,r_+]=\{W+\e <  - C \lambda^{-\frac23} \e^{\frac13}\},
\]
we use an elliptic estimate. Denote 
\[
\omega = |W+\epsilon| + |\epsilon|^{\frac13} \lambda^{-\frac23}.
\]
Then multiplying the equation
\eqref{pseqe} by $- \lambda u$ and integrating by parts we obtain
\begin{equation}
\begin{split}
\int_{r_-}^{r_+}  \lambda  |\partial_\rs u|^2 
+ \lambda^3 |W + \epsilon| | u|^2 
d\rs =& \  \int_{r_-}^{r_+} - \lambda u g d\rs +  \lambda u
  u_{\rs} \Big|_{r_-}^{r^+}
\\
\lesssim &\ \lambda \| \omega^{\frac14} u\|_{L^\infty(r_-,r_+)}
  \|\omega^{-\frac14}g\|_{L^1} + 
\| \omega^{-\frac14}g\|_{L^1}^2 ,
\end{split} \label{k6}
\end{equation}
where for the boundary terms at $r_{\pm}$ we have used the previously
obtained pointwise bounds.  On the other hand from the fundamental
theorem of calculus one obtains
\[
\lambda^2\| \omega^{\frac14} u\|_{L^\infty(r_-,r_+)}^2 \lesssim
\int_{r_-}^{r_+} \lambda |\partial_\rs u|^2 + \lambda^3 |W + \epsilon|
| u|^2 d\rs ,
\]
where the bound \eqref{wwrs} is used for the derivative of $W$ in $[r_-,r_+]$.
Combining the last two inequalities gives the 
desired bound for $u$,
\begin{equation}
\lambda \| \omega^{\frac14} u\|_{L^\infty(r_-,r_+)} \lesssim 
 \|\omega^{-\frac14}g\|_{L^1}.
\label{k7} \end{equation}
Returning to \eqref{k6}, it also follows that 
\begin{equation}
\int_{r_-}^{r_+}  \lambda  |\partial_\rs u|^2 
+ \lambda^3 |W + \epsilon| | u|^2 
d\rs \lesssim \| \omega^{-\frac14}g\|_{L^1}^2
\label{k8} \end{equation}

It remains to obtain the pointwise bound for $u_\rs$. In 
$[r_-,r_+]$  we have $W < |\epsilon|$ therefore
$W_{\rs} \lesssim |\epsilon|^\frac12$. Given $r_0^* \in 
[r_-,r_+]$ we consider an interval $r_0^* \in I \subset [r_-,r_+]$
of size $|I| = c\lambda^{-1} \omega(r^*_0)^{-\frac12}$ with a small $c$. In $I$ the size of the weight $\omega$ is constant; indeed, $\omega$ can change at most by 
\[
|I| |W_{\rs}| = c \lambda^{-1} \omega(r^*_0)^{-\frac12} \epsilon^\frac12
\lesssim c\omega(r_0^*)
\]
where at the last step we have used the bound
$\omega(r_0^*) \geq |\epsilon|^{\frac13} \lambda^{-\frac23}$.

Within $I$ we first use the $L^2$ bound \eqref{k8} to 
estimate the average $u_\rs^I$ of $u_\rs$ in $I$,
\[
|u_\rs^I|^2 \lesssim |I|^{-1} \int_I |u_\rs|^2 d\rs 
\lesssim \omega(r_0^*)^{\frac12} 
\| \omega^{-\frac14}g\|_{L^1}^2 
\]

It remains to compute the variation of $u_{\rs}$ in $I$,
which is estimated using the equation $\partial_{\rs}^2 + \lambda^2(W+\epsilon) = g$ and \eqref{k7},
\[
\begin{split}
 \int_{I} |\partial_\rs^2 u| d\rs \lesssim  \int_{I} \lambda^2 \omega |u|  
+ |g| d\rs
\lesssim \omega(r_0^*)^{\frac14} \|\omega^{-\frac14} g\|_{L^1} 
\end{split}
\]
Together, the last two bounds show that
\[
 |u_{\rs}(r^*_0)|\lesssim \omega(r^*_0)^\frac14\|\omega^{-\frac14} g\|_{L^1}.
\]
The proof of the lemma is concluded.
\end{proof}

We continue with part (b) of the proposition.  We switch to the
Regge-Wheeler coordinates. By taking a spherical harmonics expansion
it suffices to prove the result at a fixed spherical frequency
$\lambda$. Let $g_\lambda$ be at spherical frequency $\lambda$ with
support in $\{ 11 M/4 < r < 4M\}$. Using a time frequency multiplier
with smooth symbol we can split $g_\lambda$ into two components, one
with high ($\gg \lambda$) time frequency and one with low time frequency.
We consider the two cases separately.

{\bf Case I}. $g_\lambda$ is localized at time frequencies $\{ |\tau|
\gg (1+\lambda) \}$.  This corresponds to Cases 1,2,3 in the proof of
part (a).  As a consequence of the results there we have the a-priori
bound
\[
(|\tau| + \lambda) \| u\|_{L^2} \lesssim \|(\partial_{\rs}^2 +
V_{\lambda,\tau} ) u\|_{L^2}  
\]
for all $u$ with support in $\{ 5M/2 < r < 5M\}$. By duality this
implies that for each $g \in L^2$ with support in $\{ 5M/2 < r < 5M\}$
there exists a solution $v$ to 
\[
(\partial_{\rs}^2 + V_{\lambda,\tau} ) v = g
\]
 in\footnote{no boundary condition is imposed on $v$}
 $\{ 5M/2 < r < 5M\}$ with 
\[
(|\tau| + \lambda) \| v\|_{L^2} \lesssim \|g\|_{L^2}.
\]
Applying this at all time frequencies $ |\tau| \gg (1+\lambda)$
we find a solution  $u_\lambda$ to 
\begin{equation}
L_{RW} u_\lambda = g_\lambda
\label{lrweasy}\end{equation}
 in $\{ 5M/2 < r < 5M\}$ so that
\[
(1+\lambda) \|u_\lambda \|_{L^2} +\|\partial_t u_\lambda\|_{L^2} 
\lesssim \|g_\lambda\|_{L^2} .
\]
Multiplying equation \eqref{lrweasy} 
by $\chi^2_{ps} u_\lambda$ and integrating by parts
we obtain
\[
 \|\partial_\rs (\chi_{ps} u_\lambda)\|_{L^2}^2  \lesssim 
 \lambda^2 
\|\chi_{ps} u_\lambda\|^2_{L^2}    + \|\chi_{ps} \partial_t
u_\lambda\|_{L^2}^2 + \|u_\lambda\|_{L^2}^2 + \|g_\lambda\|_{L^2}^2.
\]
Hence the function $v_\lambda = \chi_{ps} u_\lambda$ satisfies
\begin{equation}
\| \nabla v_\lambda\|_{L^2}   \lesssim \|g_\lambda \|_{L^2}.
\label{wes}\end{equation}
On the other hand, since $g_\lambda$ is supported in the smaller
interval $\{ 11 M/4 < r < 4M\}$, it follows that $v_\lambda$ solves
the equation
\[
L_{RW} v_\lambda - g_\lambda = 
[L_{RW},\chi_{ps}] u_\lambda.
\]
Here the right hand side is supported in a region, away from the
photon sphere, where the $L^2$ and $LE_{ps}^*$ norms are equivalent.
Then this is seen to satisfy
\[
\| L_{RW} v_\lambda - g_\lambda\|_{L^2} \lesssim
\|g_\lambda\|_{L^2}
\]
by applying \eqref{wes} with $\chi_{ps}$ replaced by a cutoff with
slightly larger support.

Finally, the standard energy estimates for $v_\lambda$ allow us to
obtain uniform energy bounds for $v_\lambda$ from the averaged energy
bounds in \eqref{wes}, thus improving \eqref{wes} to
\begin{equation}
  \| \nabla v_\lambda\|_{L^2} +  \| \nabla v_\lambda\|_{L^\infty L^2}
 \lesssim \|g_\lambda \|_{L^2}.
\label{wesa}\end{equation}

{\bf Case II}. $g_\lambda$ is localized at time frequencies $\{ |\tau|
\lesssim (1+\lambda) \}$.  This corresponds to Case 4 in the proof of part
(a). We first observe that the result in part (a) can be strengthened
to
\begin{equation}
 \| \phi\|_{LE_{ps}}^2\lesssim \| f\|^2_{ LE_{ps}^* + L^1 L^2}.
\label{psar}\end{equation}
Indeed, suppose that $f = f_1+ f_2$ with $f_1 \in  LE_{ps}^*$ 
and $f_2 \in L^1 L^2$. We  
 solve the forward problem
\[
\Box_g \phi_2 = f_2.
\]
By Theorem~\ref{theorem.1} and Duhamel's formula we have
\[
\|\phi_2\|_{LE_0} \lesssim \|f_2\|_{L^1 L^2}.
\]
We truncate $\phi_2 \to \tilde \chi_{ps}(r) \phi_2$ in a slightly
larger set than the support of $f_2$ and compute
\[
\| \Box_g  (\tilde \chi_{ps}\phi_2) - f_2\|_{LE_{ps}^*} = 
\| [\Box_g,  \tilde \chi_{ps}]\phi_2\|_{LE_{ps}^*}
\approx \|  [\Box_g,  \tilde \chi_{ps}]\phi_2\|_{L^2}
\lesssim \|\phi_2\|_{LE_0}
\]
since the above commutator is supported in a compact set in $r$ 
away from the photon sphere.

From Duhamel's formula and part (a) of the proposition it follows that
\[
\|  \tilde \chi_{ps}\phi_2 \|_{LE_{ps}} \lesssim \|f_2\|_{L^1 L^2}.
\]
On the other hand applying directly part (a) of the proposition 
to $\phi - \tilde \chi_{ps}\phi_2$ we obtain
\[
\|  \phi - \tilde \chi_{ps}\phi_2 \|_{LE_{ps}} \lesssim 
\| \Box_{g} ( \phi - \tilde \chi_{ps}\phi_2 )\|_{LE_{ps}^*}
\lesssim \|f_1\|_{LE_{ps}^*}+
\|f_2\|_{L^1 L^2}.
\]
Hence \eqref{psar} follows.

 As a consequence of \eqref{psar} we obtain
\[
\lambda \| (a_{ps}^{-1})^w(\lambda)  u_\lambda\|_{L^2}
\lesssim \inf_{L_{RW} u_\lambda = g_1+ g_2} \Bigl( 
\| a_{ps}^w(\lambda) g_1\|_{L^2} 
+ \|g_2\|_{L^1 L^2}\Bigr).
\]
By duality, from this bound from below for $L_{RW}$, we obtain 
a local solvability result. Precisely, for each $g_\lambda$
at spherical frequency $\lambda$ with support 
in $\{ 5M/2 < r < 5M\}$  there is a function
$u_\lambda$ in the same set which solves
\begin{equation}
L_{RW} u_\lambda = g_\lambda
\label{frecleq}\end{equation}
and satisfies the bound
\begin{equation}
\lambda( \| (a_{ps}^{-1})^w(\lambda)  u_\lambda\|_{L^2}
+ \| u_\lambda\|_{L^\infty L^2})  \lesssim 
\| a_{ps}^w(\lambda) g_{\lambda}\|_{L^2}.
\label{frecl}\end{equation}
Since $(a_{ps}^{-1})^w$ has an inverse in $OPS^\delta_{1,0}$, from the
first term above we also obtain an $L^2$ bound for $u_\lambda$, namely
\begin{equation}
\lambda^{1-\delta}\| u_\lambda\|_{L^2} \lesssim 
\| a_{ps}^w(\lambda) g_{\lambda}\|_{L^2}.
\label{freclerror}\end{equation}
Since  $g_\lambda$ is localized at time frequencies $|\tau| \lesssim
(1+\lambda)$, it follows that $u_\lambda$ above can be assumed 
to have a similar time frequency localization. Hence \eqref{frecl}
also gives
\begin{equation}
  \| (a_{ps}^{-1})^w(\lambda)
  u_{\lambda t}\|_{L^2}
+ \|  u_{\lambda t}\|_{L^\infty L^2} \lesssim 
\| a_{ps}^w(\lambda) g_{\lambda}\|_{L^2}. 
\label{frecla}\end{equation}
We can also obtain a similar bound for the $\rs$ derivative of
$u_\lambda$.   For the local energy part we multiply \eqref{frecleq} by
$\chi_{ps}((a_{ps}^{-1})^w(\lambda))^2 \chi_{ps} u_\lambda$. After some
commutations where all errors are bounded using the previous estimates we
obtain
\[
\begin{split}
\|(a_{ps}^{-1})^w(\lambda) \partial_\rs (\chi_{ps} u_\lambda)\|_{L^2}^2 
\! \lesssim &\  \lambda^2
\|(a_{ps}^{-1})^w(\lambda) \chi_{ps} u_\lambda\|_{L^2}^2 +
 \|(a_{ps}^{-1})^w (\lambda)\chi_{ps} u_{\lambda t} \|_{L^2}^2  \\ &\  + 
\lambda^{2-2\delta}  \|u_\lambda\|^2_{L^2}  + \| g_\lambda \|_{L^2}^2. 
\end{split}
\]
For the $L^\infty L^2$ bound on $ \partial_{\rs}(\chi_{ps} u_\lambda)$ we consider a
smooth compactly supported function $\chi(t)$. Then multiplying
\eqref{frecleq} by $\chi^2 \chi_{ps}^2 u_\lambda$ and commuting we
obtain 
\[
\| \chi \partial_\rs (\chi_{ps} u_\lambda)\|_{L^2}^2 
\lesssim  \lambda^2
\|\chi \chi_{ps} u_\lambda\|_{L^2}^2 +
 \| \chi  \chi_{ps} u_{\lambda t} \|_{L^2}^2    + 
  \|u_\lambda\|^2_{L^2}  + \| g_\lambda \|_{L^2}^2.
\]
Taking also \eqref{frecl} and \eqref{frecla} into account we have a
bound on local averaged energy for $ \chi \chi_{ps} u_\lambda$:
\[
\|  \partial_\rs  ( \chi \chi_{ps} u_\lambda)\|_{L^2}^2 
+  \lambda^2\|\chi \chi_{ps} u_\lambda\|_{L^2}^2 +
 \| \partial_t (\chi  \chi_{ps} u_{\lambda}) \|_{L^2}^2    \lesssim 
  \| a_{ps}^w(\lambda)  g_\lambda \|_{L^2}^2.
\]
By energy estimates applied to $ \chi \chi_{ps} u_\lambda$ we can
convert the averaged energy bound into a pointwise energy bound to
obtain
\[
\|  \partial_\rs  ( \chi \chi_{ps} u_\lambda)\|_{L^\infty L^2}^2 
+  \lambda^2\|\chi \chi_{ps} u_\lambda\|_{L^\infty L^2}^2 +
 \| \partial_t (\chi  \chi_{ps} u_{\lambda}) \|_{L^\infty L^2}^2    \lesssim 
  \| a_{ps}^w(\lambda)  g_\lambda \|_{L^2}^2.
\]
Summing up \eqref{frecl}, \eqref{frecla} and the similar bounds above
for the $\rs$ derivatives we finally obtain
\[
\|(a_{ps}^{-1})^w(\lambda) \nabla ( \chi_{ps} u_\lambda)\|_{L^2} + \| \nabla (
\chi_{ps} u_\lambda)\|_{L^\infty L^2} \lesssim \|g_\lambda\|_{LE^*} ,
\]
where $\nabla = (\partial_\rs, \partial_t,\lambda)$.

On the other hand if $g_\lambda$ is supported in $\{11M/4 < r < 4M\}$
then $u_\lambda$ solves the equation
\[
L_{RW} \chi_{ps} u_\lambda  - g_\lambda = [L_{RW}, \chi_{ps}]  u_\lambda .
\]
The right hand side is supported away from the photon sphere, 
where the $L^2$ and $LE_{0}^*$ norms are equivalent.
Then, by applying Theorem \ref{theorem.1} with $\chi_{ps}$ replaced by a cutoff with
slightly larger support, this is seen to satisfy 
\[
\|L_{RW} \chi_{ps} u_\lambda - g_\lambda\|_{LE_0^*} \lesssim
\|g_\lambda\|_{LE_{ps}^*}.
\]
The proof of the proposition is concluded.
\end{proof}

\subsection{ The analysis at infinity}
\label{inftysec}

In the Schwarzschild space $\M$, if a function
$u$ in $\M$ is supported in $\{ r > 4M\}$ we interpret it as a
function in $\R \times \R^{3}$ by setting $u(t, x) = u(t,r,\omega)$
for $x = r \omega$.   We now state the analogue of
Proposition~\ref{pps}.

\begin{proposition}
a) Let $\phi$ solve $\Box_g \phi = 0$  in $\{r > 4M\}$. Then
\[
\| \chi_{\infty} \phi\|_{LE_M}^2 \lesssim \| \phi\|_{LE_0}^2 + E[\phi](0).
\]

b) Let $f \in LE^*_M$ be supported in $\{r > 4M\}$. Then
there  is a function $\phi$ supported in $\{ r > 3M\}$
which solves  $\Box_g \phi = f$ in $\{r \gg M\}$
so that
\begin{equation}
\sup_t E[\phi](t) + \| \phi\|_{LE_M}^2 + 
\| \Box_g \phi - f\|_{L^2}^2  \lesssim \| f\|^2_{ LE^*_M}.
\label{infsolve}\end{equation}

\label{pinfty}\end{proposition}

\begin{proof}
a) For $R > 0$ we denote by $\chi_{>R}$ a smooth cutoff function 
which is supported in $\{|x| > R\}$ and equals $1$ in $\{|x| \geq 2R\}$. 
If $R > 4M$ then 
\[
\|  (\chi_{\infty} - \chi_{> R}) \phi\|_{LE_{M}}^2 \lesssim \| \phi\|_{LE_0}^2.
\]
It remains to show that for a fixed sufficiently large $R$ we have
\[
\| \chi_{> R} \phi\|_{LE_M}^2 \lesssim \| \phi\|_{LE_0}^2 + E[\phi](0).
\]
For this we notice that $\chi_{>R} \phi$ solves the equation
\begin{equation}\label{InftyEqn} 
\Box_g( \chi_{> R}\phi) = f_1(x) \nabla \phi + f_2(x)\phi ,
\end{equation}
where $f_1$ and $f_2$ are supported in $\{ R < |x| < 2R\}$.  If $R$ is
sufficiently large then outside the ball $\{ |x| \leq R\}$ the
operator $\Box_g$ is a small long range perturbation of the
d'Alembertian.  Then the estimate \eqref{Mest} applies, see
e.g. \cite[Proposition 2.2]{MT1} or \cite[(2.23)]{MS2} (with no
obstacle, $\Omega = \emptyset$) and we have
\[
\begin{split}
\| \chi_{> R} \phi\|_{LE_M}^2 \lesssim  &\  E[\chi_{>R} \phi](0)
+ \| \Box_g( \chi_{>R}\phi)\|^2_{LE^*_M}
\\
\lesssim  &\  E[ \phi](0)
+ \| [\Box_g, \chi_{>R}]\phi\|^2_{L^2}
\\
\lesssim  &\  E[ \phi](0)
+ \| \phi \|^2_{LE_0} ,
\end{split}
\]
where in the last two steps we have used the compact support of $
\Box_g( \chi_{>R}\phi) = [\Box_g, \chi_{>R}]\phi$.

b) Let $R$ be large enough, as in part (a). For $|x| > R$ the
Schwarzchild metric $g$ is a small long range perturbation of the
Minkowski metric, according to the definition in \cite{MT}.  We
consider a second metric $\tilde{g}$ in $\R^{3+1}$ which coincides
with $g$ in $\{|x| > R\}$ but which is globally a small long range
perturbation of the Minkowski metric. Let $\psi$ be the forward solution 
to $\Box_{\tilde g} \psi = f$. Then we set 
\[
\phi = \chi_{>R} \psi.
\]
 The estimate \eqref{Mest} holds for the metric $\tilde{g}$, therefore
we obtain
\[
\sup_t E[\psi](t) + \| \psi\|_{LE_M} \lesssim \| f\|^2_{ LE^*_M}
\]
Then the same bound holds as well for $\phi$. Furthermore, we can
compute the error 
\[
\Box_g \phi - f = (\chi_{>R} -1) f +  [\Box_g, \chi_{>R}] \psi
\]
This has compact spatial support, and can be easily estimated in $L^2$ 
as in part (a).

\end{proof}

\subsection{ Proof of Theorem~\ref{theorem.3}}
\label{tproof}
Given $f \in LE^*$ we split it into
\[
f = \chi_{eh} f + \chi_{ps} f + \chi_\infty f.
\]
For  the last two terms we use part (b) of the
Propositions~\ref{pps},\ref{pinfty} to produce
approximate solutions   $\phi_{ps}$ and $\phi_\infty$ 
near the photon sphere, respectively near infinity.
Adding them up we obtain an approximate solution 
\[
\phi_0 = \phi_{ps} + \phi_\infty
\]
for the equation $\Box_g \phi = f$.  Due to \eqref{psb} and
\eqref{infsolve} we obtain for $\phi_0$ the bound
\begin{equation}
\sup_\tv E[\phi_0](\tv) + 
\|\phi_0 \|_{LE}^2 
\lesssim \| f\|^2_{ LE^*},
\label{infsolvea}\end{equation}
while the error 
\[
f_1 = \Box_g (\phi _{ps} + \phi_\infty)- f
\]
is supported away from $r = 3M$ and $r = \infty$ and satisfies
\[
\|f_1\|_{ LE^*_0} \approx \|f_1\|_{L^2} \lesssim \|f\|_{ LE^*}.
\]

Then we find $\phi = \phi_0+\phi_1$ by solving
\[
\Box_g \phi_1 = f_1 \in  LE^*_0, \qquad \phi_1[0] = \phi[0]-\phi_0[0],
\]
By Theorem~\ref{theorem.1} we obtain the $LE_0$ bound for $\phi_1$.
It remains to improve this to an $LE$ bound for $\phi_1$.  By part (a)
of Proposition~\ref{pinfty} we can estimate $\| \chi_{\infty}
\phi\|_{LE_M}$. 

Near the photon sphere we would like to apply part (a) of
Proposition~\ref{pps} to $\chi_{ps} \phi$. However we cannot proceed
in an identical manner because part (a) of Proposition~\ref{pps} does
not involve the Cauchy data of $\phi$ at $t = 0$, and instead applies
to functions $\phi$ defined on the full real axis in $t$. To address this issue 
we extend $\phi_1$ backward in $t$ to the set $\M'_R$, by solving the 
homogeneous problem $\Box_g \phi_1 = 0$ in $\M'_R$, with matching 
Cauchy data on the common boundary of $\M_R$ and $\M_R'$.
The extended function $\phi_1$ belongs to both $LE(\M_R)$ and $LE(M_R')$,
and now we can estimate $\chi_{ps} \phi_1$ 
via part (a) of Proposition~\ref{pps}.

\section{Strichartz estimates}

In this section we prove Theorem~\ref{Strichartz.theorem}.
The theorem follows from the following two propositions. The first
gives the result for the right hand side, $f$, in the dual local energy space:

\begin{proposition} \label{pl2tolp} Let $(\rho,p,q)$ be a nonsharp
  Strichartz pair. Then for each $\phi \in LE$ with $\Box_g \phi \in
  LE^* + L^1_\tv L^2$ we have
\begin{equation}
\|\nabla \phi\|_{L^{p}_\tv \dot H^{-\rho,q}}^2
\lesssim  E[\phi](0)+ \|\phi\|_{LE}^2+ \|\Box_g \phi\|_{LE^*+L^1_\tv L^2}^2.
\end{equation}
\end{proposition}

The second one allows us to use $L^{p_2'} L^{q_2'}$ in the right hand
side of the wave equation.

\begin{proposition} \label{plqtolp}
There is a parametrix $K$ for $\Box_g$ so that for all  nonsharp
  Strichartz pairs $(\rho_1,p_1,q_1)$ and $(\rho_2,p_2,q_2)$ we have
\begin{equation}
\sup_{\tv} E[Kf](\tv) + E[Kf](\Sigma_R^+) + \| Kf\|_{LE}^2 + \| \nabla Kf
\|_{L^{p_1}_\tv \dot H^{-\rho_1,q_1}}^2 \lesssim \| f\|_{L^{p_2'}_{\tv}
  \dot H^{\rho_2,q_2'}}^2
\label{kfest}\end{equation}
and the error estimate
\begin{equation}
  \| \Box_g Kf - f\|_{LE^*+L^1_\tv L^2}
 \lesssim \| f\|_{L^{p_2'}_\tv \dot H^{\rho_2,q_2'}}.
\label{kferr}\end{equation}
\end{proposition}

We first show how to use the propositions in order to prove the Theorem.

\begin{proof}[Proof of Theorem~\ref{Strichartz.theorem}]
Suppose that $\Box_g \phi = f$ with $f \in L^{p_2'} \dot H^{\rho_2,q_2'}$.
We write $\phi$ as 
\[
\phi = \phi_1 + Kf
\]
with $K$ as in Proposition~\ref{plqtolp}. By \eqref{kfest}
the $Kf$ term satisfies all the required estimates; therefore it
remains to consider $\phi_1$. Using also \eqref{kferr} we obtain
\[
\|\Box_g \phi_1\|_{LE^*+L^1_\tv L^2}^2 + E[\phi_1](0) \lesssim E[\phi](0) +
\|f\|_{L^{p_2'} \dot H^{\rho_2,q_2'}}^2.
\]
Then Theorem~\ref{theorem.3} combined with Duhamel's formula yields
\[
\|\phi_1\|_{LE}^2 +\|\Box_g \phi_1\|_{LE^*+L^1_\tv L^2}^2 
+\sup_{\tilde
      v} E[\phi_1](\tilde v) \lesssim
E[\phi](0) + \|f\|_{L^{p_2'}\dot H^{\rho_2,q_2'}}^2.
\]
Finally the $L^{p_1} \dot H^{-\rho_1,q_1}$ bound for $\nabla \phi_1$ follows by 
Proposition~\ref{pl2tolp}.\end{proof}

We continue with the proofs of the two propositions.

\begin{proof}[Proof of Proposition~\ref{pl2tolp}]
By Duhamel's formula and Theorem~\ref{theorem.3} we can neglect 
the $L^1 L^2$ part of $\Box_g \phi$. Hence in the sequel we assume 
that $\Box_g \phi \in LE^*$.

We use cutoffs to split the space into three regions, namely near the 
event horizon, near the photon sphere and near infinity,
\[
\phi = \chi_{eh} \phi + \chi_{ps} \phi + \chi_{\infty} \phi.
\]
Due to the definition of the $LE$ and $LE^*$ norms we have
\[
\begin{split}
  E[\phi](0) + \|\phi\|_{LE}^2 + \|\Box_g \phi\|_{LE^*}^2 \gtrsim &\
  E[ \chi_{eh}\phi](0) + \| \chi_{eh}\phi\|_{H^1}^2 + \|\Box_g (
  \chi_{eh}\phi)\|_{L^2}^2
  \\
  + &\ E[ \chi_{ps}\phi](0) + \| \chi_{ps}\phi\|_{LE_{ps}}^2 + \|\Box_g (
  \chi_{ps}\phi)\|_{LE^*_{ps}}^2
 \\
+ &\ E[ \chi_{\infty}\phi](0) + \| \chi_{\infty}\phi\|_{LE_M}^2 + \|\Box_g (
  \chi_{\infty}\phi)\|_{LE^*_M}^2.
\end{split}
\]
Proving this requires commuting $\Box_g$ with the cutoffs. However
this is straightforward since the $LE$ and $LE^*$ norms are equivalent
to the $H^1$, respectively $L^2$, norm in the support of $\nabla
\chi_{eh}$, $\nabla \chi_{ps}$ and $\nabla \chi_{\infty}$.

It remains to prove the $L^p_{\tv} \dot H^{-\rho,q}$ bound for each of the 
three terms in $\nabla \phi$. We consider the three cases separately:

{\bf I. The estimate near the event horizon.} This is the easiest
case. Given $\phi$ supported in $\{ r < 11M/4\}$, we partition it on
the unit scale with respect to $\tv$,
\[
\phi = \sum_{j \in \Z} \chi(\tv-j) \phi ,
\]
where $\chi$ is a suitable smooth compactly supported bump function.
Commuting the cutoffs with $\Box_g$ one easily obtains the square
summability relation
\[
\sum_{j \in \N} \| \chi(\tv-j) \phi\|_{H^1}^2 + \| \Box_g (\chi(\tv-j)
\phi)\|_{L^2}^2  + E[\chi(\tv-j)
\phi](0) \lesssim \|\phi\|_{H^1}^2  + \| \Box_g \phi\|_{L^2}^2 + E[\phi](0), 
\]
where the energy term on the left is nonzero only for finitely many $j$.
Since each of the functions $\chi(\tv-j) \phi$ have compact support, 
they satisfy the Strichartz estimates due to the local theory;
see \cite{MSS}, \cite{Smith}, \cite{T3}. The above square summability 
with respect to $j$ guarantees that the local estimates can be added
up. 

{\bf II. The estimate near the photon sphere.}
For $\phi$ supported in $\{ 5M/2 < r < 5M\}$ we need to show that
\[
\|\nabla \phi\|_{L^p_\tv H^{-\rho,q}}^2 \lesssim
E[\phi](0)+\|\phi\|_{LE_{ps}}^2 + \|\Box_g \phi\|_{LE^*_{ps}}^2.
\]
We use again the Regge-Wheeler coordinates. Then the operator 
$\Box_g$ is replaced by $L_{RW}$. The potential $V$ can
be neglected due to the straightforward bound
\[
\|V \phi\|_{LE^*_{ps}} \lesssim \|\phi\|_{LE_{ps}}.
\]
Indeed, for $\phi$ at spherical frequency $\lambda$ we have
\[
\|V \phi\|_{LE^*_{ps}} \lesssim |\ln (2+\lambda)|^\frac12 \|\phi\|_{L^2}
\lesssim \lambda |\ln (2+\lambda)|^{-\frac12} \|\phi\|_{L^2} \lesssim
\|\phi\|_{LE_{ps}}.
\]

 We introduce the auxiliary function
\[
\psi = A_{ps}^{-1} \phi.
\]
By the definition of the $LE_{ps}$ norm we have
\begin{equation}
\| \psi \|_{H^1} \lesssim \|\phi\|_{LE_{ps}}.
\label{h1psi}\end{equation}
We also claim that 
\begin{equation}
\|L_{RW} \psi\|_{L^2}  \lesssim \| \phi\|_{LE_{ps}} + \| L_{RW} \phi\|_{L^2}.
\label{ppsi}\end{equation}
Since $A_{ps}^{-1}$ is $L^2$ bounded, this is a consequence of the 
commutator bound
\[
[  A_{ps}^{-1}, L_{RW}]: LE_{ps} \to L^2,
\]
or equivalently
\begin{equation}
[  A_{ps}^{-1}, L_{RW}] A_{ps}: H^1 \to L^2.
\label{cmt}\end{equation}
It suffices to consider the first term in the symbol calculus, as the
remainder belongs to $OPS^\delta_{1,\delta}$, mapping $H^\delta$ to
$L^2$ for all $\delta > 0$. The symbol of the first term is 
\[
q(\xi,\rs,\lambda) = \{ a_{ps}^{-1}(\lambda), \xi^2 + r^{-3}(r-2M) \lambda^2\}
a_{ps}(\lambda)
\]
and a-priori we have $q \in S^{1+\delta}_{1,\delta}$.
For a better estimate we compute the Poisson bracket
\[
q(\xi,\rs,\lambda)  = a_{ps}^{-1}(\lambda) \gamma_y(y,\ln \lambda) 
\frac{ 4 \xi \rs - 2 \xi \partial_{\rs}( r^{-3}(r-2M))}{\rs^2 +
  \lambda^{-2} \xi^2}
\]
where $y = \rs^2 + \lambda^{-2} \xi^2$. The first two factors on the
right are bounded. The third is bounded by $\lambda$ since
$\partial_{\rs}( r^{-3}(r-2M))$ vanishes at $\rs = 0$. In addition,
$q$ is supported in $|\xi| \lesssim \lambda$.  Hence we obtain $q \in
\lambda S^{0}_{1-\delta,\delta}$. Then the commutator bound
\eqref{cmt} follows.

Given \eqref{h1psi} and \eqref{ppsi}, we argue as in the first case,
namely we localize $\psi$ to time intervals of unit length and then
apply the local Strichartz estimates. By summing over these strips we
obtain
\[
\| \nabla \psi \|_{L^p H^{-\rho,q}} \lesssim \| \phi\|_{LE_{ps}} + \|
L_{RW} \phi\|_{L^2}
\]
for all sharp Strichartz pairs $(\rho,p,q)$.

To  return to $\phi$ we invert $A_{ps}^{-1}$,
\[
\phi = A_{ps} \psi + (1- A_{ps} A_{ps}^{-1}) \phi.
\]
The second term is much more regular,
\[
\| \nabla (1- A_{ps} A_{ps}^{-1}) \phi\|_{L^2 H^{1-\delta}} \lesssim
\| \phi \|_{LE_{ps}}, \qquad \delta > 0;
\]
therefore it satisfies all the Strichartz estimates simply by Sobolev 
 embeddings. 

For the main term $A_{ps}\psi$ we take advantage of the fact that
we only seek to prove the nonsharp Strichartz estimates for $\phi$.
The  nonsharp Strichartz estimates for $\psi$ are obtained from the
sharp ones via Sobolev embeddings,
\[
\|\nabla \psi\|_{H^{-\rho_2,q_2}} \lesssim \|\nabla
\psi\|_{H^{-\rho_1,q_1}}, \qquad \frac{3}q_2 + \rho_2 = \frac{3}q_1 +
\rho_1, \quad \rho_1 < \rho_2.
\]
To obtain the nonsharp estimates for $\phi$ instead, we need a
slightly stronger form of the above bound, namely

\begin{lemma}
Assume that $1 < q_1 < q_2 < \infty$. Then
\begin{equation}\label{ps.Sobolev}
\|A_{ps} u \|_{H^{-\rho_2,q_2}} \lesssim \| u\|_{H^{-\rho_1,q_1}}, \qquad 
\frac{3}q_2 + \rho_2 = \frac{3}q_1 + \rho_1.
\end{equation}
\end{lemma}
\begin{proof}
We need to prove that the operator
\[
B = Op^w(\xi^2+\lambda^2+1)^{-\frac{\rho_2}2} A_{ps} 
Op^w(\xi^2+\lambda^2+1)^{\frac{\rho_1}2}
\] 
maps $L^{q_1}$ into $L^{q_2}$. The principal symbol of 
$B$ is 
\[
b_0(\rs,\xi,\lambda) = (\xi^2+\lambda^2+1)^{\frac{\rho_1-\rho_2}2}
a_{ps}(\rs,\xi,\lambda),
\]
and by the pdo calculus the remainder is easy to estimate,
\[
B - b_0^w \in OPS^{ \rho_1-\rho_2 -1 +\delta}_{1,0}, \qquad \delta > 0.
\]
The conclusion of the lemma will follow from the
Hardy-Littlewood-Sobolev inequality if we prove a suitable pointwise
bound on the kernel $K$ of $b_0^w$, namely
\begin{equation}
|K(r^*_1,\omega_1,r^*_2,\omega_2)| \lesssim (|r^*_1-r^*_2| 
|\omega_1-\omega_2|^2)^{-1 + \frac{1}{q_1} -\frac{1}{q_2}}.
\label{Kbd}\end{equation}
For fixed $\rs$ we consider a smooth dyadic partition of unity 
in frequency as follows:
\[
1  = \chi_{\{|\xi| > \lambda\}} 
+ \sum_{\mu\, \text{dyadic}} \chi_{\{\lambda \approx \mu\}}
\left(\chi_{\{ |\xi| \lesssim \nu_0\}} 
+ \sum_{\nu = \nu_0}^\mu \chi_{\{
    |\xi| \approx \nu\}}\right) ,
\]
 where $\nu_0 = \nu_0(\lambda,\rs)$ is given by
\[
\ln \nu_0(\lambda,\rs) = \ln \lambda +  \max\{ \ln \rs, -\sqrt{ \ln{\lambda}}\}.
\]
This leads to a similar decomposition for $b_0$, namely
\[
b_0 = b_{00} + \sum_\mu \left( b_{\mu, < \nu_0} + \sum_{\nu = \nu_0}^\mu
b_{\mu \nu} \right).
\]
In the region $|\xi| \gtrsim \lambda$ the symbol $b_0$ is of class
$S^{\rho_1-\rho_2}$, which yields a kernel bound  for $b_{00}$
of the form 
\[
|K_{00}(r^*_1,\omega_1,r^*_2,\omega_2)| \lesssim (|r^*_1-r^*_2| +
|\omega_1-\omega_2|)^{-3 -\rho_1+\rho_2}.
\]
The symbols of $b_{\mu \nu}$ are supported in $\{ |\xi| \approx \nu, \
\lambda \approx \mu \}$, are smooth on the same scale and have
size $\ln (\nu^{-1} \mu) \mu^{\rho_2-\rho_1}$. 
Hence their kernels satisfy bounds 
of the form
\[
|K_{\mu,\nu} (r^*_1,\omega_1,r^*_2,\omega_2)| \lesssim \ln (\nu^{-1}
\mu)  \mu^{\rho_1-\rho_2} \nu (|r^*_1-r^*_2|\nu + 1)^{-N} \mu^2
(|\omega_1-\omega_2|\mu + 1)^{-N}
\]
and similarly for $K_{\mu,<\nu_0}$. Then \eqref{Kbd} follows after
summation.
\end{proof}

{\bf III. The estimate near infinity.}

Let us first recall the setup from \cite{MT}. We fix a Littlewood-Paley dyadic decomposition of frequency space in $\R^3$,
\[ 
1=\sum_{k=-\infty}^\infty S_k(D),\quad \text{supp}\, s_k\subset
\{2^{k-1}<|\xi|<2^{k+1}\}.
\]
Functions $u$ in $\R \times \R^3$ which are localized to frequency
$2^k$ are measured in
\begin{equation}\label{X_knorm}
\|u\|_{X_k} = 2^{k/2}\|u\|_{L^2(A_{<-k})} + \sup_{j\ge -k}
\||x|^{-1/2} u\|_{L^2(A_j)} ,
\end{equation}
where
\[ 
A_j = \R\times \{2^j\le |x|\le 2^{j+1}\},\quad A_{<j}=\R\times
\{|x|\le 2^j\}.
\]
As  in \cite{MT}, by $X^0$ we denote the space of functions in $\R\times \R^3$
with norm
\begin{equation}\label{X0norm}
\| u\|_{X^0}^2 = \sum_{k=-\infty}^\infty \|S_k u\|^2_{X_k}
\end{equation}
and by $Y^0$ the dual norm
\[
\| u\|_{Y^0}^2 = \sum_{k=-\infty}^{\infty} \|S_k u\|^2_{X_k^{'}}
\]
where $X_k^{'}$ is the dual norm of $X_k$.

 One can establish the following (see \cite[Lemma 1]{MT})
\begin{lemma}\label{LEMvsX^0}
 The following inequalities hold:
\begin{equation}\label{LEM<X^0}
\sup_j 2^{-j/2} \|\nabla u\|_{L^2(A_j)}
\lesssim \|\nabla u\|_{X^0}
\end{equation}
and its dual
\begin{equation}\label{LE_M^*>Y^0}
\|u \|_{Y^0} \lesssim \|u \|_{LE^*_M}.
\end{equation}
\end{lemma}

For small deviations from the Minkowski metric, one can also establish
stronger local energy estimates involving the $X^0$ and $Y^0$ norms;
more precisely, one can prove (see \cite[Theorem 4]{MT}):
\begin{lemma}\label{X^0LE}
  Let $\tilde{g}$ be a sufficiently small, long range perturbation of
  the Minkowski metric.  Then, for all solutions $u$ to the
  inhomogeneous problem $\Box_{\tilde{g}} u = f$ one has
\[
 \|\nabla u\|_{L^{\infty}_t L^2_x}^2 + \| \nabla u\|_{X^0}^2
 \lesssim E[u](0) + \|f\|_{Y^0+L^1_t L^2_x}^2.
\] 
 \end{lemma}

We now return to proving our estimate. For $\phi$ supported in $\{ r > 4M\}$ we need to show that
\begin{equation}
\|\nabla \phi\|_{L^p \dot H^{-\rho,q}}^2 \lesssim
E[\phi](0)+\|\phi\|_{LE_M}^2 + \|\Box_g \phi\|_{LE^*_M}^2.
\label{phin}\end{equation}
For large $R$ we split $\phi$ into a near and a far part
\[
\phi = \chi_{>R} \phi + \chi_{<R} \phi
\]
and estimate
\[
\begin{split}
  E[\phi](0) + \|\phi\|_{LE_M}^2 + \|\Box_g \phi\|_{LE^*_M}^2 \gtrsim &\
  E[ \chi_{>R}\phi](0) + \| \chi_{>R}\phi\|_{LE_M}^2 + 
\|\Box_g (\chi_{>R}\phi)\|_{LE^*_M}^2
  \\
  + &\ E[ \chi_{<R}\phi](0) + \| \chi_{<R}\phi\|_{H^1}^2 + \|\Box_g (
  \chi_{<R}\phi)\|_{L^2}^2.
\end{split}
\]
The term $\chi_{<R} \phi$ has compact support in $r$ and can be
treated as in the first case (i.e. near the event horizon). Hence
without any restriction in generality we can restrict ourselves to the
case when $\phi$ is supported in $\{ r > R\}$. But in this region the
operator $\Box_g$ is a small long range perturbation of $\Box$;
therefore the results of \cite{MT} apply. More
precisely, 
from \cite[Theorem 7(a)]{MT} we obtain
\[
 \|\nabla \phi\|_{L^p \dot H^{-\rho,q}}^2 \lesssim
E[\phi](0)+\|\nabla \phi\|_{X^0}^2 + \|\Box_g \phi\|_{LE^*}^2.
\] 
 This does not directly imply \eqref{phin}, since the $X^0$ norm is stronger than $LE_M$. However, we can apply Lemma \ref{X^0LE} and \eqref{LE_M^*>Y^0} to obtain the bound
\[
 \|\nabla \phi\|_{X^0}^2 \lesssim E[\phi](0) + \|\Box_g \phi\|_{LE^*_M}^2.
\]

\end{proof}

\begin{proof}[Proof of Proposition~\ref{plqtolp}]
We split $f$ into 
\[
f = \chi_{eh} f + \chi_{ps} f +\chi_\infty f
\]
and construct the parametrix separately in the three regions.

{\bf I. The parametrix near the event horizon.}
 We further partition the term $\chi_{eh} f$ into unit intervals
\[
\chi_{eh} f = \sum_j \chi(\tv-j) \chi_{eh} f
\]
with $\chi$ supported in $[-1,1]$, so that each component has compact
support in the region 
\[
D_j = \{ r_0 \leq r < 11M/4,\ j-2 < \tv < j+2 \}.
\]
Let $\psi_j$ be the forward solution to
\[
\Box_g \psi_j = \chi(\tv-j) \chi_{eh} f.
\]
Due to the local Strichartz estimates for variable coefficient wave
equations, we obtain the uniform bounds
\[
\| \nabla \psi_j\|_{L^{p_1} H^{-\rho_1,q_1}(D_j) } +\|\nabla
\psi_j\|_{L^\infty L^2(D_j)} + \|\psi_j\|_{L^\infty L^2(D_j)} \lesssim 
\|  \chi(\tv-j) \chi_{eh} f\|_{L^{p'_2} H^{\rho_2,q'_2}}.
\]
Next we truncate $\psi_j$ using a cutoff function $\tilde \chi(\tv
-j,r)$ which is supported in $D_j$ and equals $1$ in the support 
of $\chi(\tv-j) \chi_{eh}$. Then the bound above also holds for the
truncated functions $\phi_j = \tilde \chi(\tv
-j,r) \psi_j$,
\begin{equation}
\| \nabla \phi_j\|_{L^{p_1} H^{-\rho_1,q_1} } +\|\nabla
\phi_j\|_{L^\infty L^2} +  \|\phi_j\|_{L^\infty L^2(D_j)}  \lesssim 
\|  \chi(\tv-j) \chi_{eh} f\|_{L^{p'_2} H^{\rho_2,q'_2}}.
\label{ehp1} \end{equation}
In addition,
\[
\Box_g \phi_j - \chi(\tv-j) \chi_{eh} f = [\Box_g,  \tilde \chi(\tv
-j,r)] \psi_j;
\]
therefore
\begin{equation}
\| \Box_g \phi_j - \chi(\tv-j) \chi_{eh} f\|_{L^2}  \lesssim 
\|  \chi(\tv-j) \chi_{eh} f\|_{L^{p'_2} H^{\rho_2,q'_2}}.
\label{ehp2} \end{equation}
Finally, by energy estimates we also obtain a bound for the energy of 
$\phi_j$ on the future space-like boundary of $D_j$ at $r=r_0$,
\begin{equation}
\|\nabla \phi_j\|_{L^2(D_j \cap \{r=r_0\})} + \|\phi_j\|_{L^2 (D_j \cap \{r=r_0\})}\lesssim 
\|  \chi(\tv-j) \chi_{eh} f\|_{L^{p'_2} H^{\rho_2,q'_2}}.
\label{ehp3} \end{equation}

To conclude we set 
\[
K_{eh} f = \sum_j \phi_j.
\]
Summing up the bounds \eqref{ehp1}, \eqref{ehp2} and \eqref{ehp3}
for $\phi_j$ we obtain the desired bounds for $K_{eh}$, namely
\[
\sup_\tv E[K_{eh} f](\tv) + E[K_{eh} f](\Sigma_R^+) +\| K_{eh}
f\|_{H^1}^2+ \| \nabla K_{eh} f\|_{L^{p_1} H^{-\rho_1,q_1}}^2 \lesssim 
\| \chi_{eh} f\|_{L^{p'_2} H^{\rho_2,q'_2}}^2,
\]
respectively the error estimate
\[
\| \Box_g K_{eh} f -\chi_{eh} f\|_{L^2}  \lesssim 
\| \chi_{eh} f\|_{L^{p'_2} H^{\rho_2,q'_2}}.
\]

{\bf II. The parametrix near the photon sphere.}
We work in the Regge-Wheeler coordinates. Arguing as in the previous
case we produce a parametrix $\tilde K_{ps}$ with the property that,
for each $f$ supported in $\{ 5M/2 +\epsilon < r < 5M-\epsilon\}$, the
function $\tilde K_{ps} f$ is supported in $\{ 5M/2 < r < 5M\}$
 and satisfies the bounds
\[
\sup_t E[\tilde K_{ps} f](t)  +\| \tilde K_{ps}
f\|_{H^1_{x,t}}^2+ \| \nabla \tilde K_{ps} f\|_{L^{p_1} H^{-\rho_1,q_1}}^2 \lesssim 
\|  f\|_{L^{p'_2} H^{\rho_2,q'_2}}^2
\]
and the error estimate
\[
\| L_{RW}  \tilde K_{ps} f - f\|_{L^2}  \lesssim 
\| f\|_{L^{p'_2} H^{\rho_2,q'_2}}.
\]
Then we define the localized parametrix near the photon sphere
$K_{ps}$ as 
\[
K_{ps} f = A_{ps}^{-1} \tilde K_{ps} \tchi_{ps} A_{ps} (\chi_{ps} f)
\]
with $\tchi_{ps} =1$ in the support of $\chi_{ps}$ and slightly larger
support.  Then we  show that $K_{ps}$ satisfies the required bounds.

 We recall that $(\rho_2,p_2,q_2)$ is a nonsharp  Strichartz pair.
Then by  \eqref{ps.Sobolev} we can write
\[
\| \tchi_{ps}A_{ps} (\chi_{ps} f)\|_{L^{p'_3} H^{\rho_3,q'_3}}
\lesssim \|\chi_{ps} f\|_{L^{p'_2} H^{\rho_2,q'_2}}
\]
for some other Strichartz pair $(\rho_3,p_3,q_3)$ with $p_3=p_2$ and
$q_3 < q_2$.
Since $A_{ps}^{-1}$ is $L^2$ bounded, from the above bounds for $\tilde K_{ps}$ 
we obtain
\[
\sup_t E[K_{ps} f](t)  +\|  K_{ps}
f\|_{H^1}^2 \lesssim \|\chi_{ps} f\|_{L^{p'_2} H^{\rho_2,q'_2}}^2.
\]
By using \eqref{ps.Sobolev} with $A_{ps}$ replaced by the weaker 
operator $A_{ps}^{-1}$ we also obtain the $L^{p_1} H^{-\rho_1,q_1}$
bound for $K_{ps} f$:
\[
 \| \nabla  K_{ps} f\|_{L^{p_1} H^{-\rho_1,q_1}}^2 \lesssim 
\| \nabla \tilde K_{ps} \tchi_{ps}  A_{ps}(\chi_{ps} f)\|_{L^{p} H^{-\rho,q}}^2 
\lesssim   \|\chi_{ps} f\|_{L^{p'_2} H^{\rho_2,q'_2}}^2 ,
\]
where $(\rho,p,q)$ is another  Strichartz pair with $p=p_1$ and $q < q_1$.

It remains to consider the error estimate, 
\begin{equation}
\| L_{RW} K_{ps} f - \chi_{ps} f\|_{LE^{*}+L^1_\tv L^2} \lesssim \|\chi_{ps}
f\|_{L^{p'_2} H^{\rho_2,q'_2}} ,
\label{erlrw}\end{equation}
for which we compute
\[
\begin{split}
  L_{RW} K_{ps} f - \chi_{ps} f =&\ [L_{RW}, A_{ps}^{-1}] \tilde K_{ps}
  \tchi_{ps} A_{ps}(\chi_{ps} f)
\\
&\ + A_{ps}^{-1} ( L_{RW} \tilde K_{ps} - I) \tchi_{ps} A_{ps} (\chi_{ps} f)
\\
&\ + (A_{ps}^{-1} \tchi_{ps} A_{ps} - \tchi_{ps}) (\chi_{ps} f).
\end{split}
\]
We consider each term in the above decomposition. For the first term,
due to the $H^1$ bound for $\tilde{K}$, we need the commutator bound
\[
[L_{RW},A_{ps}^{-1}] : H^1 \to LE^*
\]
or equivalently
\[
A_{ps} [L_{RW},A_{ps}^{-1}] : H^1 \to L^2 ,
\]
which is almost identical to \eqref{cmt} and is proved in the same manner.

The bound for the second term is a direct consequence of the 
$L^2$ error bound for $\tilde K$.

Finally, for the last term we know that
$(A_{ps}^{-1} A_{ps} - I) \in OPS^{-1+\delta}_{1,0}$; therefore using
Sobolev embeddings we estimate
\[
\| (A_{ps}^{-1}\tchi_{ps} A_{ps} - \tchi_{ps}) (\chi_{ps} f)\|_{L^{p'_2} H^\frac12}
  \lesssim   \|\chi_{ps}f\|_{L^{p'_2} H^{\rho_2,q'_2}}.
\]
This concludes the proof of \eqref{erlrw} since
\[
L^{p'_2} H^\frac12 \subset L^{2} H^\frac12+L^{1} H^\frac12
\subset LE^*_{ps} + L^1 L^2.
\]

{\bf III. The parametrix near infinity.}

We now consider the last component of $f$, namely $\chi_\infty f$.
For some large $R$ we separate it into two parts,
\[
\chi_\infty f =( \chi_\infty-\chi_{>R})  f + \chi_{>R} f.
\]
The first part has compact support in $r$; therefore we can handle it
as in the first case (i.e. near the event horizon), producing a
parametrix $K_\infty^{<R}$. For the second part we modify the metric
$g$ for $r < R$ to a metric $\tilde g$ which is a small, long-range
perturbation of $\Box$. We let $\psi_{\infty}$ be the forward solution
to
\[
\Box_{\tilde g}\psi_{\infty} = \chi_{>R} f.
\]
We consider a second cutoff function $\tilde \chi_{>R}$ which is
supported in $r > R$ and equals $1$ in the support of $\chi_{>R}$.
Then we define
\[
K_\infty^{>R} f = \tilde \chi_{>R} \psi_\infty .
\]
It remains to show that  $K_\infty^{>R} $ satisfies the
appropriate bounds, 
\[
\sup_t E[K_\infty^{>R} f](t) + \| K_\infty^{>R} f \|^2_{LE_M} +
\|\nabla K_\infty^{>R} f \|_{L^{p_1} \dot H^{-\rho_1,q_1}}^2 \lesssim
\|\chi_{>R} f \|_{L^{p'_2} \dot H^{\rho_2,q'_2}}^2,
\]
respectively the error estimate
\[
\| \Box_g K_\infty^{>R} f - \chi_{>R}f\|_{LE^*_M} \lesssim \|\chi_{>R} f
\|_{L^{p'_2} \dot H^{\rho_2,q'_2}}.
\]
These are easily obtained by applying the following lemma to $\psi_\infty$:

\begin{lemma}
Let $f \in L^{p'_2} \dot H^{\rho_2,q'_2}$. Then the forward solution 
$\psi$ to $\Box_{\tilde g} \psi = f$ satisfies the  bound
\begin{equation}
\sup_t E[ \psi](t)  + \| \psi\|_{LE_M}^2 + \|\nabla \psi \|_{L^{p_1}
 \dot H^{-\rho_1,q_1}}^2 \lesssim \| f\|_{L^{p'_2} \dot H^{\rho_2,q'_2}}^2.
\label{t6mt}\end{equation}
\end{lemma}

It remains to prove the lemma. This largely follows from  
\cite[Theorem 6]{MT}, but there is an interesting technical issue
that needs clarification. Precisely, \cite[Theorem 6]{MT} shows that 
we have the bound
\begin{equation}
\sup_t E[ \psi](t)  + \| \nabla \psi\|_{X^0}^2 + \|\nabla \psi \|_{L^{p_1}
 \dot H^{-\rho_1,q_1}}^2 \lesssim \| f\|_{L^{p'_2} \dot H^{\rho_2,q'_2}}^2.
\label{backref}\end{equation}

By Lemma \ref{LEMvsX^0}, we are left with proving that
\begin{equation}\label{ll2} 
\sup_{j\in \Z} 2^{-\frac{3j}2}\|\psi\|_{L^2 (A_j)} \lesssim \| f\|_{L^{p_2'} \dot H^{\rho_2,q'_2}}.
\end{equation}

 We note that this does not follows from Lemma \ref{LEMvsX^0}; this is a forbidden endpoint of the Hardy
inequality in \cite[Lemma 1(b)]{MT}.


However, the bound \eqref{ll2} can still be obtained, although 
in a roundabout way. Precisely, from \eqref{backref} we have
\begin{equation}
\sup_t E[ \psi](t) \lesssim \| f\|_{L^{p_2'} \dot H^{\rho_2,q'_2}}^2
\label{f1}\end{equation}
for the forward in time evolution, and similarly for the backward in time
problem.

On the other hand, a straightforward modification of the classical
Morawetz estimates (see e.g. \cite{MS}) for the wave equation shows that
the solutions to the homogeneous wave equation $\Box_{\tilde g} \psi = 0$
satisfy 
\begin{equation}
\sup_{j\in \Z} 2^{-3j}\|\psi\|_{L^2 (A_j)}^2 \lesssim  E[ \psi](0).
\label{f2}\end{equation}
Denote by $1_{t>s} H(t,s)$ the forward fundamental solution for $\Box_{\tilde g}$
and by $H(t,s)$ its backward extension to a solution to the
homogeneous equation, $\Box_{\tilde g} H(t,s) = 0$.
Combining the bounds \eqref{f1} and \eqref{f2} shows that
\[
\sup_j 2^{-3j}\left\| \int_\R H(t,s) f(s) ds\right\|_{L^2 (A_j)}^2  \lesssim 
\| f\|_{L^{p'_2} \dot H^{\rho_2,q'_2}}^2.
\]
Since $p'_2 < 2$, by the Christ-Kiselev lemma \cite{CK}, it follows that
\[
\sup_j 2^{-3j} \left\| \int_t^\infty  H(t,s) f(s) ds \right\|_{L^2 (A_j)}^2  \lesssim 
\| f\|_{L^{p'_2} \dot H^{\rho_2,q'_2}}^2 ,
\]
which is exactly \eqref{ll2}.

\end{proof}


\section{ The critical NLW}

In this section we prove Theorem~\ref{tnlw}.  We first consider
\eqref{nonlin} in the compact region $\M_C$. We denote by $\psi$ the
solution to the homogeneous equation
\[
\Box_g \psi = 0, 
\qquad \psi_{|\Sigma_0} = \phi_0,\ \tilde K \psi_{|\Sigma_0} = \phi_1
\]
and by $Tf$ the solution to the inhomogeneous problem
\[
\Box_g (Tf) = f, \qquad 
TF_{|\Sigma_0} = 0,\ \tilde K Tf_{|\Sigma_0} = 0.
\]
Then we can rewrite the nonlinear equation \eqref{nonlin}
in the form
\begin{equation}
\phi = \psi \pm  T(\phi^5).
\label{fixp}\end{equation}

We define Sobolev spaces in $\M_C$ by restricting to $\M_C$ functions
in the same Sobolev space which are compactly supported in a larger
open set.  By the local Strichartz estimates we have
\[
\| \psi\|_{H^{\frac12,4}(\M_C)} \lesssim E[\phi](\Sigma_0) 
\]
and
\[
\|Tf\|_{H^{\frac12,4}(\M_C)} \lesssim \|f\|_{H^{\frac12,\frac43}(\M_C)}.
\]
At the same time we have the multiplicative estimate
\[
\| \phi^5 \|_{H^{\frac12,\frac43}(\M_C)} \lesssim
\|\phi\|_{H^{\frac12,4}(\M_C)}^5.
\]
Then for small initial data we can use the contraction principle to solve 
\eqref{fixp} and obtain a solution $\phi \in H^{\frac12,4}(\M_C)$.
In addition, still by local Strichartz estimates, the solution $\phi$
will have finite energy on any space-like surface, in particular
on the forward and backward space-like boundary of $\M_C$.
Thus we obtain 
\[
E[\phi](\Sigma_R^-) \lesssim E[\phi](\Sigma_0).
\]

It remains to solve \eqref{nonlin} in $\M_R$ (and its other three
symmetrical copies). Using the $(\tv,r,\omega)$ coordinates in 
$\M_R$ we define $\psi$ and $T$ as above, but with 
Cauchy data on $\Sigma_R^-$.

By the global Strichartz estimates in Theorem~\ref{Strichartz.theorem},
for $(s,p)$  as in the theorem we have 
\[
\| \psi\|_{L^p \dot H^{s,p}(\M_R)} \lesssim E[\phi](\Sigma_R^-) 
\]
and 
\[
\| Tf \|_{L^p \dot H^{s,p}(\M_R)} \lesssim \|f\|_{L^1 L^2}.
\]

In particular we can take  $p=5$ which corresponds to $s =
\frac{3}{10}$. By Sobolev embeddings we have
\[
\| \phi\|_{L^5 L^{10}} \lesssim \|\phi\|_{\dot H^{\frac3{10},5}};
\]
therefore 
\[
\| \phi^5\|_{L^1 L^{2}} \lesssim \|\phi\|_{\dot H^{\frac3{10},5}}^5.
\]
Hence we can solve \eqref{fixp} using the contraction 
principle and obtain a solution $\phi \in {\dot H^{\frac3{10},5}}$.
This implies that $\phi^5 \in L^1 L^2$, which yields all of the 
other Strichartz estimates, as well as the energy bound 
on the forward boundary $\Sigma_R^+$ of $\M_R$.
This concludes the proof of the theorem.

\bigskip

\end{document}